\input ./style/arxiv-general.cfg
\documentclass[aos,MSNbibl,seceqn,dvips]{arximspdf}
\makeatletter
   \@ifpackageloaded{graphicx}{}{\usepackage{graphicx}}
\makeatother
\doi{10.1214/15-AOS1328}% Updated by VTEXPTS2LaTeX.exe, 24.03.2015
\volume{43}
\issue{5}
\pubyear{2015}
\firstpage{1865}
\lastpage{1895}
\docsubty{FLA}

\makeatletter
\def\cal{\mathcal}
\newtheorem{thmm}{Theorem}
\newproclaim{rmk}{Remark}
\newtheorem{lem}{Lemma}
\newproclaim{exa}{Example}
\newtheorem{cor}{Corollary}
\newcommand{\wtd}{\widetilde}
\newcommand{\cN}{{\cal N}}
\makeatother

\begin{document}
\begin{frontmatter}

\title{Optimal detection of multi-sample aligned sparse~signals}
\runtitle{Optimal detection aligned sparse signals}

\begin{aug}
% Corresponding author: Hock Peng Chan - stachp@nus.edu.sg% Updated by
%VTEXPTS2LaTeX.exe, 24.03.2015 12:55
\author[A]{\fnms{Hock Peng}~\snm{Chan}\corref{}\ead[label=e1]{stachp@nus.edu.sg}\thanksref{T1}}
\and
\author[B]{\fnms{Guenther}~\snm{Walther}\ead[label=e2]{gwalther@stat.stanford.edu}\thanksref{T2}}
\runauthor{H.~P. Chan and G. Walther}
\thankstext{T1}{Suported by the National University of Singapore Grant
R-155-000-120-112.}
\address[A]{Department of Statistics\\
\quad and Applied Probability\\
National University of Singapore\\
6 Science Drive 2 \\
Singapore 117546 \\
\printead{e1}}
\affiliation{National University of Singapore and Stanford University}
\thankstext{T2}{Supported by NSF Grants DMS-10-07722 and DMS-12-20311.}
\address[B]{Statistics Department\\
Stanford University\\
390 Serra Mall\\
Stanford, California 94305\\
USA\\
\printead{e2}}
%\author[A]{\fnms{}~\snm{}\corref{}\ead[label=e1]{}}%,
%\author[]{\fnms{}~\snm{}\ead[label=]{}}
% \and
%\author[]{\fnms{}~\snm{}\ead[label=]{}}
%\runauthor{}
%\affiliation{}
%\dedicated{}
%\address[A]{\\\printead{e1}}
%\address[]{\\\printead{}}
\end{aug}

% HISTORY:
%
\received{\smonth{12} \syear{2014}}% Updated by VTEXPTS2LaTeX.exe,
%24.03.2015 12:55
%
\revised{\smonth{2} \syear{2015}}% Updated by VTEXPTS2LaTeX.exe,
%24.03.2015 12:55

% ABSTRACT
%
\begin{abstract}
We describe, in the detection of multi-sample aligned sparse signals,
the critical boundary separating detectable from nondetectable signals,
and construct tests that achieve optimal detectability:
penalized versions of the Berk--Jones and the higher-criticism test
statistics evaluated over pooled scans,
and an average likelihood ratio over the critical boundary.
We show in our results an inter-play between the scale of the sequence
length to signal length ratio,
and the sparseness of the signals.
In particular the difficulty of the detection problem is not noticeably
affected unless this ratio grows exponentially with the number of sequences.
We also recover the multiscale and sparse mixture testing problems as
illustrative special cases.
\end{abstract}

% KEYWORDS

%
%\begin{keyword}[class=AMS]
%\kwd[Primary ]{}
%\kwd{}
%\kwd[; secondary ]{}
%\end{keyword}
%\begin{keyword}
%\kwd{}
%\end{keyword}
\end{frontmatter}

%s1 #&#
\section{Introduction}\label{sec1}
Consider a population of sequences having a common time (or location) index.
Signals, when they occur,
are present in a small fraction of the sequences and aligned in time.
In the detection of copy number variants (CNV) in multiple DNA sequences,
Efron and Zhang \cite{EZ10} used local f.d.r.,
Zhang et al. \cite{ZSJL10} and Siegmund, Yakir and Zhang \cite{SYZ11}
applied scans of weighted $\chi^2$-statistics,
Jeng, Cai and Li \cite{JCL13} applied higher-criticism test statistics.
Tartakovsky and Veeravalli \cite{TV08}, Mei \cite{Mei10} and Xie and
Siegmund \cite{XS13}
considered the analogous sequential detection of sparse aligned changes
of distribution in parallel streams of data,
with applications in communications, disease surveillance, engineering
and hospital management.
These advances have brought in an added multi-sample dimension
to traditional scan statistics works (see, e.g., the papers in \cite{GPW08})
that consider a single stream of data.
%to
%traditional scan statistics works;
%see, for example, \cite{GPW08}
%for the consideration of a single stream of data.

In this paper, we tackle the problem of detectability of aligned sparse signals,
extending sparse mixture detection (cf. \cite{BJ79,CW12,DJ04,HJ10,Ing97,Ing98,Wal13}) to aligned signals,
and extending multiscale detection (cf. \cite{CW13,DS01,LT00,Roh09})
to multiple sequences.
Hence not surprisingly,
we incorporate ideas developed by the sparse mixture and multiscale
detection communities
to find the critical boundary separating detectable from nondetectable
hypotheses.\vadjust{\goodbreak}
In Arias-Castro, Donoho and Huo \cite{ADH06a,ADH06},
there are also links between sparse mixtures and multiscale detection methods
in the detection of a sparse component on an unknown low-dimensional
curve within a higher-dimensional space.
Our work here is less geometrical in nature
as the aligned-signal assumption allows us to reduce the problem to one
dimension by summarizing across sample first.

We supply optimal adaptive max-type tests:
penalized scans of the higher criticism and Berk--Jones test statistics.
We also supply an optimal Bayesian test:
an average likelihood ratio (ALR) that tests against alternatives lying
on the critical boundary.
The rationale behind the ALR is to focus testing at the most sensitive
parameter values,
where small perturbations can result in sharp differences of detection powers.

We state the main results in Section~\ref{sec2}.
We describe the detectable region of aligned sparse signals in the
multi-sample setting,
and show that the penalized scans achieve asymptotic detection power 1 there.
We learn from the detection boundary the surprising result that the
requirement to locate the signal in the time domain
does not affect the overall difficulty of the detection problem,
unless the sequence length to signal length ratio grows exponentially
with the number of sequences.

In Section~\ref{sec3}, we show the optimality of the ALR and consider special cases
of our model that have been well studied in the literature using
max-type tests:
the detection of a signal with unknown location and scale in a single sequence,
and the detection of a sparse mixture in many sequences of length 1.
We show that the general form of our ALR provides optimal detection in
these important special cases.
We also illustrate the detectability and detection
of multi-sample signals on a CNV dataset.

In Section~\ref{extensions}, the detection problem is extended to heteroscedastic signals.
The extension illustrates the adaptivity of the penalized scans.
Even though the detection boundary has to be extended to take into
account the heteroscedasticity,
the penalized scans as described in Section~\ref{sec2} are still optimal.
On the other hand, the ALR tests have to be re-designed to ensure
optimality under heteroscedasticity.
The model set-up here is similar to that in Jeng, Cai and Li \cite{JCL13}.
There optimality is possible without imposing penalties on the scan of
the higher-criticism test
because the signal length was assumed to be very short.

%s2 #&#
\section{Main results}\label{sec2}
Let $\{ (X_{n1}, \ldots, X_{nT})\dvtx 1 \leq n \leq N \}$ be a population
of sequences.
We consider the prototypical set-up
%
%e2.1 #&#
\begin{equation}
\label{proto} X_{nt} = \mu_{nt}+Z_{nt} \qquad
\mbox{where } Z_{nt} \mbox{ are i.i.d. } \mathrm{N}(0,1).
\end{equation}
Under the
null hypothesis $H_0$ of no signals, $\mu_{nt}=0$ for all $n$ and
$t$. Under the alternative hypothesis $H_1$ of aligned signals,
there exists an unknown $q>0$ of disjoint intervals $(j_T^{(k)},
j_T^{(k)}+\ell_T^{(k)}]$ such that for the $k$th interval,
$1 \leq k \leq q$,
there is a probability $\pi_N^{(k)}>0$ that this interval has an
elevated mean
%
%e2.2 #&#
\begin{eqnarray}
\label{munt} \mu_{nt} & = & \cases{ %
\mu_N^{(k)} I_n^{(k)}/\sqrt{
\ell_T^{(k)}}, & \quad $\mbox{if } j_T^{(k)}
< t \leq j_T^{(k)}+\ell_T^{(k)},$
\vspace*{2pt}\cr
0, & \quad$\mbox{otherwise,}$}
\nonumber
\\[-8pt]
\\[-8pt]
\nonumber
I_n^{(k)} & \sim& \operatorname{ Bernoulli} \bigl(
\pi_N^{(k)} \bigr),
\end{eqnarray}
with $\mu_N^{(k)} >0$ and the $I_n^{(k)}$'s and $Z_{nt}$'s jointly independent.
Let $\pi_N = \pi_N^{(1)}$, $\mu_N=\mu_N^{(1)}$ and so forth.

Model (\ref{munt}) extends sparse mixture detection by adding a
time-dimension, and there is a similar extension in potential
applications. For example, in the
detection of bioweapons use, as introduced
in \cite{DJ04}, we can assume that
there are $N$ observational units in a geographical region, each
accumulating information over time on
bioweapons usage. The bioweapons are in use over a specific but unknown
time period, and only a small fraction of the units are
affected. Alternatively in covert communications detection, only a
small fraction of $N$ detectors, each tuned to a distinct signal
spectrum, observes
unusual activities during the period in which communications are taking
place. In the detection of genes
that are linked to cancer, readings of DNA
copy numbers are taken from the chromosomes of $N$ cancer patients, and
only a small fraction of the patients exhibits
copy number changes at the gene locations.
In Section~\ref{extensions}, we shall consider an extension of (\ref{proto}) and
(\ref
{munt}) to signals carrying a noise component.

In the detection of copy number changes, the common practice was to
process samples one at a time; see Lai et al. \cite{Lai05}.
In contrast, Efron and Zhang \cite{EZ10}, Zhang et al. \cite{ZSJL10}
and Jeng et al. \cite{JCL13}
proposed procedures that pool across samples first. Our analysis here
shows that the alignment information is
important, and we should indeed pool across samples first.
In Appendix \ref{appC}, we provide a
comparison between pooling information across sample versus pooling
information within sample first.

Consider $\pi_N = N^{-\beta}$ for some $\frac{1}{2} < \beta< 1$.
Ingster \cite{Ing97,Ing98} and Donoho and Jin
\cite{DJ04} showed that in the special case $T=1$ (hence $q=1$,
$j_1=0$, $\ell_1=1$),
as $N \rightarrow\infty$, the
critical detectable value of $\mu_N$ is $b_N^*(\beta)
:= \sqrt{2 \rho^*(\beta) \log N}$, where
%
%e2.3 #&#
\begin{equation}
\label{rhos} \rho^*(\beta) = \cases{ %
\beta-
\tfrac{1}{2}, & \quad$\mbox{if } \tfrac{1}{2} < \beta\leq\tfrac{3}{4},$
\vspace*{2pt}\cr
(1-\sqrt{1-\beta})^2, & \quad$\mbox{if } \frac{3}{4} < \beta< 1.$}
\end{equation}
That is, if $\mu_N= \sqrt{2 \rho\log N}$ with $\rho< \rho^*$,
then no test can detect that $\mu_N \neq0$ in the sense that the sum
of Type I and Type II error probabilities tends to 1 for any test.
Donoho and Jin \cite{DJ04} further showed that Tukey's higher criticism
as well as the Berk--Jones statistic
achieve the detection boundary $b_N^*$;
that is, if $\rho> \rho^*$, then the sum of Types I and~II error
probabilities tends to 0.
Jager and Wellner \cite{JW07} showed that Tukey's test is a member of a
family of goodness-of-fit $\phi$-divergence tests
that can each achieve the detection boundary.

When $T>1$, we need to deal with the complication of multiple comparisons
over $j_T$ and $\ell_T$, and the question arises of how much harder the
detection problem becomes. The
number of disjoint intervals in $(0,T]$ with common length
$\ell_T$ is approximately $T/\ell_T$. This ratio has to be factored into
the computation of the
detection boundary. The main message of Section~\ref{sec2.1} is that
the difficulty of the detection problem is not noticeably affected unless
this ratio of sequence length to signal length grows exponentially with
the number $N$ of sequences. Sections \ref{sec2.2} and \ref{sec2.3} provide optimal
max-type tests that attain the detection boundary.

%s2.1 #&#
\subsection{Detectability of aligned signals}\label{sec2.1}

Let $a_m \sim b_m$ if $\lim_{m \rightarrow\infty} (a_m/b_m)=1$ and
$a_m \,\dot{\sim}\,  b_m$ if $\lim_{m \rightarrow\infty} (a_m/b_m)=C$
for some
constant $C>0$. Let $\lfloor\cdot\rfloor$ be the
greatest integer function and $\# B$ the number of elements in a set
$B$. Let $E_0(E_1)$ denote expectation under
$H_0(H_1)$.
We are interested here in the signal length $\ell_T^{(k)}$ in (\ref
{munt}) satisfying
%
%e2.4 #&#
\begin{equation}
\label{C1} T/\ell_T^{(k)} \sim\exp
\bigl(N^{\zeta^{(k)}}-1 \bigr)\qquad \mbox{for some } \zeta ^{(k)} \geq0.
\end{equation}
The case of $T$ varying sub-exponentially with $N$ will be considered
in Section~\ref{extensions}.

We shall show that under (\ref{C1}) with $N \rightarrow\infty$,
the asymptotic threshold detectable value of $\mu_N$ when
$\pi_N = N^{-\beta}$ and $\beta\in(0,1)$ is
%
%e2.5 #&#
\begin{eqnarray}\quad
\label{bNbz} & & b_N(\beta,\zeta)
\nonumber
\\[-8pt]
\\[-8pt]
\nonumber
& & \qquad = \cases{ %
 \sqrt{\log \bigl(1+N^{2 \beta-1+\zeta}
\bigr)}, &\quad  $\mbox{if } 0 \leq\zeta\leq1-4 \beta/3,$
\vspace*{2pt}\cr
(\sqrt{1-\zeta}-\sqrt{1-
\zeta-\beta}) \sqrt{2 \log N}, & \quad$\mbox{if } 1-4 \beta/3 < \zeta\leq1-\beta,$
\vspace*{2pt}\cr
\sqrt{N^{\beta+\zeta-1}}, &\quad $\mbox{if } \zeta> 1-\beta.$}
\end{eqnarray}

The first case $0 \leq\zeta\leq1-4 \beta/3$ can be further sub-divided
into: (a) $0 \leq\zeta\leq
1-2 \beta$, under which
%
%e2.6 #&#
\begin{equation}
\label{poly} b_N(\beta,\zeta) \,\dot{\sim}\,  N^{-(1-2 \beta-\zeta)/2}\qquad
\mbox{(decays polynomially with } N),
\end{equation}
and (b) $1-2 \beta< \zeta\leq1-4 \beta/3$, under which
%
%e2.7 #&#
\begin{equation}
\label{grow} b_N(\beta,\zeta) \sim\sqrt{(2 \beta-1+\zeta) \log N}\qquad
\mbox{(grows at } \sqrt{\log N} \mbox{ rate)}.
\end{equation}

Formula (\ref{bNbz}) specifies the functional form of $b_N$ as a
function of $\beta$.
Since $\beta$ appears in the exponent in (\ref{poly}) and in the
third case
of (\ref{bNbz}), $b_N$ is specified only up to multiplicative constants
in these cases.

The boundary $b_N$ is an extension of the Donoho--Ingster--Jin boundary $b_N^*$.
In the case of a sparse mixture,
$T=\ell_T=1$, and (\ref{C1}) is satisfied with
$\zeta=0$. By the second case in (\ref{bNbz}) and by
(\ref{grow}), $b_N(\beta,0) \sim b_N^*(\beta)$ when $\frac{1}{2} <
\beta
< 1$. Furthermore, $b_N(\beta,0)$ in (\ref{poly}) recovers the detection
boundary in the dense case $0 < \beta\leq\frac{1}{2}$ established by
Cai, Jeng and Jin \cite{CJJ11}.

Formula (\ref{bNbz}) likewise recovers the detection boundary for the
special case
of only one sequence. For the scaled mean $\mu_N$ in (\ref{munt}),
this boundary is known to be $\sqrt{2 \log(eT/\ell_T)}$
%and is
%attained, for example, by the penalized scan \cite{CW13}.
and is attained by the penalized scan; see, for example, \cite{CW13}.
To see how this special case is subsumed in the general setting above,
set $T \sim\exp(N-1)$ so that it suffices to consider $\zeta\in
(0,1)$ in (\ref{C1})
to parametrize the scale of the signal $\ell_T/T \in(0,1)$.
Then set $\beta=0$ so that the signal is present in each of the $N$ sequences.
Since the signals are aligned and have the same means,
by sufficiency one can equivalently consider the one sequence $S_t$ of
length $T$ obtained by summing the $X_{nt}$ over $n$.
Dividing by $\sqrt{N}$ to restore unit variance and formally plugging
$\beta=0$ into (\ref{poly}) gives a detection threshold for $\sqrt{N}
\mu_N$
of $\dot{\sim}  N^{\zeta/2} \sim\sqrt{\log(eT/\ell_T)}$. This yields
the above detection threshold for the one sequence problem apart from
the multiplicative constant $\sqrt{2}$, which can be recovered with
a more refined analysis in (\ref{bNbz}).

The general formula (\ref{bNbz})
shows how the growth coefficient and the phase transitions of the
$\sqrt{\log N}$ growth are altered
by the effect of multiple comparisons in the
location of signals. The formula also shows that in the case
$\zeta> 0$, the signal detection thresholds can grow polynomially
with $N$.

%th1 #&#
\begin{thmm} \label{thmm1}
Assume that (\ref{munt}) and (\ref{C1}) hold for $1 \leq k
\leq q$,
with $\mu_N^{(k)} = b_N(\beta^{(k)},\zeta^{(k)})$ and $\pi_N^{(k)} =
N^{-\beta^{(k)}-\varepsilon^{(k)}}$
for some $0 < \beta^{(k)} < 1$ and $\varepsilon^{(k)} > 0$. Under these
conditions, there is no test that
can achieve, at all $j_T^{(k)}$,
$1 \leq k \leq q$,
%
%e2.8 #&#
\begin{equation}
\label{typeI} P(\mbox{Type I error}) + P(\mbox{Type II error}) \rightarrow0.
\end{equation}
\end{thmm}

The simple likelihood ratio of $(X_{n1}, \ldots, X_{nT})$, for $H_0$
against (\ref{munt}), is $L_{n \ell_T j_T}(\pi_N,\mu_N)$, where
%
%e2.9 #&#
\begin{equation}
\label{Lnlj} L_{n \ell j} \bigl(\pi^*,\mu \bigr) = 1-\pi^*+\pi^* \exp
\bigl( \mu Y_{n \ell j}-\mu^2/2 \bigr),
\end{equation}
with $Y_{n \ell j} = \ell^{-1/2} \sum_{t=j+1}^{j+\ell} X_{nt}$.
The key to proving Theorem~\ref{thmm1} (details in Section~\ref{sec5}) is to show
that under the conditions of Theorem~\ref{thmm1},
%
%e2.10 #&#
\begin{equation}
\label{prod} \prod_{n=1}^N
L_{n \ell_T j_T}(\pi_N,\mu_N) = O_p(T/
\ell_T).
\end{equation}
That is, the likelihood ratio of the signal does not grow fast enough
to overcome the noise due to the
$\sim\! T/\ell_T$ independent comparisons
of length $\ell_T$. Theorem~\ref{thmm1} follows because the likelihood
ratio test is the most powerful test.

%s2.2 #&#
\subsection{Optimal detection with the penalized higher-criticism test}\label{sec2.2}

As an illustration, first consider sparse mixture detection. That is,
let $T=1$ and test
%
%e2.11 #&#
\begin{equation}
\label{sparse} X_n \stackrel{\mathrm{ i.i.d.}} {\sim} (1-
\pi_N) \mathrm{ N}(0,1)+\pi_N \mathrm{ N}(
\mu_N,1), \qquad 1 \leq n \leq N,
\end{equation}
for $H_0$: $\pi_N=0$ against $H_1$: $\pi_N > 0$ and $\mu_N > 0$. Let
$p_{(1)} \leq\cdots\leq p_{(N)}$ be
the ordered $p$-values of the $X_n$'s.

Donoho and Jin \cite{DJ04} proposed to separate $H_0$ from $H_1$ by
applying Tukey's higher-criticism test statistic
%
%e2.12 #&#
\begin{equation}
\label{HCN} \mathrm{ HC}_N := \max_{1 \leq n \leq({N}/{2})\dvtx p_{(n)} \geq N^{-1}}
\frac{n/N-p_{(n)}}{\sqrt{p_{(n)}(1-p_{(n)})/N}}.
\end{equation}
They showed that the higher-criticism
test is optimal for sparse mixture detection. Under $H_0$, $\mathrm{ HC}_N
\sim\sqrt{2 \log\log N}$; see \cite{DJ04}, Theorem~1. Under $H_1$,
the argument of $\mathrm{ HC}_N$ at some $p_{(n)}$
is asymptotically
larger than $\sqrt{2 \log\log N}$, when
$\pi_N = N^{-\beta}$ for some $\frac{1}{2} < \beta< 1$, and $\mu
_N$ lies
above the detection boundary
$b_N^*(\beta) = \sqrt{2 \rho^*(\beta) \log N}$. For $\mu_N$ lying below
the detection boundary, it is not possible to separate $H_0$ from $H_1$.
Cai et al. \cite{CJJ11} showed that optimality extends to $\beta\in
(0,\frac{1}{2})$.

We motivate the extension of the higher-criticism test to $T>1$ by first
considering a fixed, known signal on the interval $(j,j+\ell]$. By
sufficiency, testing for an aligned signal
there is the same as testing $H_0$ against $H_1$ for the sample
$Y_{1 \ell j}, \ldots, Y_{N \ell j}$.
Let $p_{(1) \ell j} \leq\cdots\leq p_{(N) \ell j}$ be the ordered
$p$-values of the sample, and let $s_{\ell T} =
\log(eT/\ell)$. We define the higher-criticism test statistic on this
interval to be
%
%e2.13 #&#
\begin{equation}
\label{HCNlj} \mathrm{ HC}_{N \ell j} := \max_{1 \leq n \leq({N}/{2})\dvtx p_{(n)} \geq
{s_{\ell T}}/{N}}
\frac{n/N-p_{(n) \ell j}}{\sqrt{p_{(n) \ell j}(1-p_{(n) \ell j})/N}}.
\end{equation}

For $\ell=T$, the constraint in (\ref{HCNlj}) becomes $p_{(n)} \geq
N^{-1}$, which agrees with the constraint
in (\ref{HCN}). As explained in \cite{DJ04}, Section~3,
the standardization of $p_{(n)}$ given in (\ref{HCNlj}) has
increasingly heavy tails as $n$ becomes smaller,
so if HC$_{N \ell j}$ is defined without constraints on $p_{(n)}$,
then it has large values frequently due to the smallest $p_{(n)}$.
For $\ell< T$, the multiple comparisons when maximizing HC$_{N \ell
j}$ over $j$
necessitates a more restrictive constraint of $p_{(n)} \geq s_{\ell T}/N$.

The term $s_{\ell T}$ appears also in the scan of the higher-criticism
test statistic
%
%e2.14 #&#
\begin{equation}
\label{PHC} \mathrm{PHC}_{NT} := \max_{(j,j+\ell) \in B_T} (
\mathrm{ HC}_{N \ell j} - \sqrt{s_{\ell T} \log s_{\ell T}}),
\end{equation}
as a penalty that increases with $T/\ell$ to counter-balance the generally
higher scores under $H_0$ for larger $T/\ell$ when maximizing HC$_{N
\ell j}$ over $j$.

We will now specify the scanning set $B_T$ in (\ref{PHC}).
In applications $T$ is
often large, so maximizing HC$_{N \ell j}$ over all $j$ and $\ell$ is
computationally expensive; the cost is $NT^2$. We
construct below an approximating set $B_T$, similar to
that in Walther~\cite{Wal10} and
Rivera and Walther \cite{RW12}, which has a computation cost of
$NT \log T$.

Construction of $B_T$: Let $d_{r,T}=\lfloor T/(r^{1/2} e^r) \rfloor+1$,
and let
%
%e2.15 #&#
\begin{equation}\qquad
\label{BrT} B_{r,T} = \bigl\{ (j,j+\ell) \in(d_{r,T} {
\mathbf Z})^2\dvtx 0 \leq j \leq T-\ell, T/e^r < \ell\leq
T/e^{r-1} \bigr\}.
\end{equation}
We define
$B_T = \bigcup_{r=1}^{r_T} B_{r,T}$, where $r_T = \lfloor\log T \rfloor$.
The specification
of $d_{r,T}$ is so that for any $(j_T, \ell_T)$,
we can find $(j_T^*,\ell_T^*) \in B_{r,T}$ for some $r$ such
that $j_T \leq j_T^* < j_T^*+\ell_T^* \leq j_T+\ell_T$ and
%
%e2.16 #&#
\begin{equation}
\label{llr} 1-\ell_T^*/\ell_T = O
\bigl(r^{-1/2} \bigr).
\end{equation}
This property plays a part in ensuring that the loss of information due
to restriction to $B_T$ is negligible.

%th2 #&#
\begin{thmm} \label{thmm2a}
Assume (\ref{munt}) and that for some $1 \leq k \leq q$,
(\ref{C1}) holds and $\mu_N^{(k)} = b_N(\beta^{(k)},\zeta
^{(k)})$, $\pi_N^{(k)} = N^{-\beta^{(k)}+\varepsilon^{(k)}}$
for some $0 < \beta^{(k)} < 1-\zeta^{(k)}$ and $0 < \varepsilon^{(k)}
\leq
\beta^{(k)}$.
Under these conditions, $P(\mbox{Type I error}) + P(\mbox{Type
II error})
\rightarrow0$ can be achieved by testing with $\mathrm{PHC}_{NT}$.
\end{thmm}

For signal identification, when applying the penalized higher-criticism
test statistics at a threshold $c$:
\begin{longlist}[(1)]
\item[(1)] Rank the pairs $(j,j+\ell) \in B_T$ in order of descending
values of HC$_{N \ell j}-\sqrt{s_{\ell T} \log s_{\ell T}}$,
and remove those pairs with values less than $c$.

\item[(2)] Starting with the highest-ranked pair and moving downward,
remove a pair from the list if its interval
overlaps with that of a higher-ranked pair still on the list by more
than a fraction $f \geq0$ of its length.
\end{longlist}

Jeng et al. \cite{JCL13}
focused on the detection of signal segments that are well separated.
Hence their signal identification procedure is restricted to $f=0$.
Zhang et al. \cite{ZSJL10} focused on both the detection of signal
segments that are well separated,
as well as the detection of overlapping or nested signal segments.
Hence their procedure allows for $f>0$.
If there are a finite number of well-separated signal segments,
then intuitively,
all the segments are identified with probability converging to 1,
in the sense that a segment with local maximum score, in a suitably
defined neighborhood of each signal segment, is identified.

%s2.3 #&#
\subsection{Optimal detection with the penalized Berk--Jones test}\label{sec2.3}

Let $K(x,t) = x \log( \frac{x}{t} ) +(1-x) \log
( \frac{1-x}{1-t} )$ if $x \geq t$ and $K(x,t)=0$ otherwise.
This is the Berk--Jones \cite{BJ79} test statistic that was first
proposed as a more powerful alternative to the
Kolmogorov--Smirnov test statistic for testing a distribution function;
see also Owen \cite{Owe95}.
Jager and Wellner \cite{JW07} showed that there is a class of test
statistics that includes the Berk--Jones
and higher-criticism test statistics as special cases that can be used
to detect sparse mixtures (\ref{sparse}) optimally.
Specifically for $T=1$,
the testing of $\pi_N > 0$ and $\mu_N > 0$ in (\ref{sparse}) can be
detected optimally by
%
%e2.17 #&#
\begin{equation}
\label{BJN} \mathrm{ BJ}_N := N \max_{1 \leq n \leq({N}/{2})\dvtx p_{(n)} < {n}/{N}}
K(n/N, p_{(n)}).
\end{equation}

Therefore, analogously to (\ref{HCNlj}),
%
%e2.18 #&#
\begin{equation}
\label{BJNlj} \mathrm{ BJ}_{N \ell j} := N \max_{1 \leq n \leq({N}/{2})\dvtx p_{(n) \ell
j} < {n}/{N}}
K( n/N, p_{(n) \ell j})
\end{equation}
can optimally detect aligned signals on the interval $(j,j+\ell]$. In
Theorem~\ref{thmm2b} below, we shall show
that analogously to (\ref{PHC}), the penalized Berk--Jones test statistic
%
%e2.19 #&#
\begin{equation}
\label{PBJ} \mathrm{PBJ}_{NT} := \max_{(j,j+\ell) \in B_T} (
\mathrm{BJ}_{N \ell j}- s_{\ell T} \log s_{\ell T})
\end{equation}
is optimal for aligned signals detection when the signal locations are unknown.

%th3 #&#
\begin{thmm} \label{thmm2b}
Assume (\ref{munt}) and that for some $1 \leq k \leq q$,
(\ref{C1}) holds and $\mu_N^{(k)} = b_N(\beta^{(k)},\zeta
^{(k)})$, $\pi_N^{(k)} = N^{-\beta^{(k)}+\varepsilon^{(k)}}$
for some $0 < \beta^{(k)} < 1$ and $0 < \varepsilon^{(k)} \leq\beta^{(k)}$.
Under these conditions,
\[
P(\mbox{Type I error}) + P(\mbox{Type II error}) \rightarrow0
\]
can be achieved by testing with $\mathrm{PBJ}_{NT}$.
\end{thmm}

As in Section~\ref{sec2.2}, a sequential approach can be used to identify
signals when the penalized Berk--Jones exceeds
a specified threshold.

%s3 #&#
\section{Optimal detection with ALR tests}\label{sec3}

We shall introduce in Section~\ref{sec3.1} an ALR that is optimal for detecting
multi-sample aligned signals. We then
consider the special cases of detecting a sparse mixture
($T=1$ with $N \rightarrow\infty$) in Section~\ref{sec3.2} and
multiscale detection in a single sequence
($N=1$ with $T \rightarrow\infty$) in Section~\ref{sec3.3}.

%s3.1 #&#
\subsection{Detecting multi-sample aligned signals}\label{sec3.1}

The ALR builds upon the likelihood ratios $L_{n \ell_T j_T}(\pi_N,\mu
_N)$ as defined in (\ref{Lnlj}), first by
substituting $\mu_N$ by
its asymptotic threshold detectable value, followed by
integrating $\pi_N = N^{-\beta}$ over $\beta$ and finally by summing
over an approximating set for $\ell_T$ and $j_T$. In view
of (\ref{C1}), let $\zeta_{\ell,NT} = \log_N [\log(T/\ell)+1]$,
and let
%
%e3.1 #&#
\begin{equation}
\label{Ln2} L_{n \ell j}(\beta) = L_{n \ell j}
\bigl(N^{-\beta}, b_N(\beta,\zeta _{\ell,NT}) \bigr);
\end{equation}
see (\ref{bNbz}) for the definition of $b_N$.

In the case of the ALR, we consider
%
%e3.2 #&#
\begin{equation}
\label{ANT} A_{NT} := \frac{6}{\pi^2} \sum
_{r=1}^{r_T \vee1} \frac{1}{r^3 e^{r+1}} \sum
_{(j,j+\ell) \in B_{r,T}} \int_0^1 \Biggl[
\prod_{n=1}^N L_{n \ell j}(\beta)
\Biggr] \,d \beta,
\end{equation}
where $r_T$ and $B_{r,T}$ are given in Section~\ref{sec2.2}.
By (\ref{BrT}), $\# B_{r,T} \leq r e^{r+1}$.

The weights in
(\ref{ANT}) are chosen for the following reason:
Since $L_{n \ell j}(\beta)$ is a likelihood
ratio for $H_0$, it has expectation 1 under $H_0$. Hence it follows
from (\ref{ANT}) that
%
%e3.3 #&#
\begin{equation}
\label{6pi} E_0(A_{NT}) = \frac{6}{\pi^2} \sum
_{r=1}^{r_T} \frac{1}{r^3 e^{r+1}} (\#
B_{r,T}) \leq\frac{6}{\pi^2} \sum_{r=1}^\infty
\frac{1}{r^2} = 1.
\end{equation}
From (\ref{6pi}), it follows that under $H_0$,\vspace*{-1pt}
$A_{NT}=O_p(1)$ uniformly over $N$ and $T$. If the aligned signals
under $H_1$ are such that
$A_{NT} \stackrel{p}{\rightarrow} \infty$ as $N \rightarrow\infty$,
then $P(\mbox{Type I error}) + P(\mbox{Type II error})
\rightarrow0$ is achieved by simply
selecting rejection thresholds going to infinity slowly enough.

The ALR (\ref{ANT}) is, by its construction, optimal when $\mu_N =
b_N(\beta,\zeta_{\ell_T,NT})$. It is not designed to be optimal at
other $\mu_N$.
However, there is really no point in being optimal at smaller $\mu_N$,
where the maximum power that can be attained
is small. At larger $\mu_N$, the ALR test has power close to 1, so
there is not much more
to be gained in being optimal there. By focusing only on the boundary
detectable values, we remove the
noise due to the consideration of unproductive likelihood ratios
associated with too large and too small
$\mu_N$.

%th4 #&#
\begin{thmm} \label{thmm2}
Assume (\ref{munt}) and that for some $1 \leq k \leq q$
(\ref{C1}) holds and $\mu_N^{(k)} = b_N(\beta^{(k)},\zeta
^{(k)})$, $\pi_N^{(k)} = N^{-\beta^{(k)}+\varepsilon^{(k)}}$
for some $0 < \beta^{(k)} < 1$ and $0 < \varepsilon^{(k)} \leq\beta
^{(k)}$. Under these conditions,
$A_{NT} \stackrel{p}{\rightarrow} \infty$. Hence $P(\mbox{Type I
error}) + P(\mbox{Type II error})
\rightarrow0$ can be
achieved by testing with $A_{NT}$.
\end{thmm}

Among the three optimal tests that we propose here,
the ALR is the most intuitive, and the
proof of its optimality is also the most straightforward. However,
its computation involves the evaluation of a nonstandard integral,
and its form is closely linked to normal errors. On the other hand,
the penalized scans involve no integrations in their computations,
and the $p$-values in their expressions are not tied to normal errors.

%s3.2 #&#
\subsection{Detecting sparse mixtures}\label{sec3.2}

This setting has been
studied in \cite{DJ04} and discussed briefly in Section~\ref{sec2}. It
corresponds to the special case $T=1$ in the above theory, and
our test statistic $A_{NT}$ simplifies as follows: $T=1$ implies
$r_T=\zeta=
\zeta_{\ell,NT}=0$, and $B_T$ contains only $j=0, \ell=1$. Hence
\[
A_{N1}=\frac{6}{\pi^2 e^2} \int_0^1
\prod_{n=1}^N \bigl( 1-N^{-\beta}
+N^{-\beta} \exp \bigl\{b_N(\beta,0) X_n -
b_N^2(\beta ,0)/2 \bigr\} \bigr) \,d \beta,
\]
where
\[
b_N(\beta,0) = \cases{ %
 \sqrt{\log
\bigl(1+N^{2\beta-1} \bigr)}, &\quad $\mbox{if $0< \beta\leq\tfrac
{3}{4}$},$
\vspace*{2pt}\cr
(1-\sqrt{1-\beta}) \sqrt{2 \log N}, &\quad $\mbox{if $\tfrac{3}{4}< \beta<1$},$}
\]
is essentially the Cai--Jeng--Jin detection boundary $b_N^*(\beta) :=
N^{-{1}/{2}+\beta}$
for $\beta\in(0,\frac{1}{2})$, and the
Donoho--Ingster--Jin detection boundary $b_N^*(\beta)=\break \sqrt{2 \rho
^*(\beta) \log N}$
for $\beta\in(\frac{1}{2},1)$.

%co1 #&#
\begin{cor} \label{thmm4}
Assume (\ref{sparse}) with $\mu_N=
b_N(\beta,0)$ and
$\pi_N=N^{-\beta+\varepsilon}$ for some $0 < \beta< 1$ and $0 <
\varepsilon
\leq\beta$.
Under these conditions, $A_{N1}
\stackrel{p}{\rightarrow} \infty$ and (\ref{typeI}) can be
achieved by testing with $A_{N1}$.
\end{cor}

%s3.3 #&#
\subsection{Signal detection in a single sequence}\label{sec3.3}

Let $N=1$ and $\pi_1=1$.
The resulting $\beta=0$ is not covered by
our general theory, but it is a boundary case, and therefore it
is of interest to see whether our general
statistic $A_{NT}$ still allows optimal detection in this important
special case. For this testing problem in a single sequence with
$T \rightarrow\infty$, it is known that the critical
detectable value of $\mu_T$ is $b_T(\ell_T)$, where $b_T(\ell) =
\sqrt{2 \log(eT/\ell)}$, and that the popular scan statistic is
suboptimal except for signals on the smallest scales; see
Chan and Walther \cite{CW13}. It is also shown there that optimal
detection can be achieved by modifying the scan with the penalty
method introduced by D\"{u}mbgen and Spokoiny \cite{DS01}, or by
employing the condensed average likelihood ratio.

Note that when analyzing a single sequence we know a priori that
$\beta=0$, and therefore it makes sense to set $\beta$ to 0 in
the definition of $A_{NT}$ rather than integrating $\beta$ over $(0,1)$.
The resulting statistic is
%
%e3.4 #&#
\begin{equation}
\label{AT} A_T := \frac{6}{\pi^2} \sum
_{r=1}^{r_T \vee1} \frac{1}{r^3 e^{r+1}} \sum
_{(j,j+\ell) \in B_{r,T}} \exp \bigl[b_T(\ell) Y_{\ell j}-b_T^2(
\ell)/2 \bigr],
\end{equation}
where $Y_{\ell j} = \ell^{-1/2}
\sum_{t=j+1}^{j+\ell} X_{1t}$.

The test statistic $A_T$ is able to achieve the
detection boundary $b_T(\ell)$ simply because it optimizes detection
power there:

%th5 #&#
\begin{thmm} \label{thmm3}
If there exist $\ell_T$ and $j_T$ such that $E_1(Y_{\ell_T j_T}) =
b_T(\ell_T)+c_T$,
with $c_T \rightarrow\infty$ as $T \rightarrow\infty$,
then $A_T \stackrel{p}{\rightarrow} \infty$ and (\ref{typeI}) can
be achieved by testing with $A_T$.
\end{thmm}

%s3.4 #&#
\subsection{An example}\label{sec3.4}

Efron and Zhang
\cite{EZ10} applied local f.d.r. to detect CNV in multi-sample DNA
sequences. Measurements from $T=42\mbox{,}075$ probes
were taken on each chromosome 1 of $N=207$ glioblastoma
subjects from the Cancer Genome Atlas Project
\cite{TCGA08}. At each probe on each sequence,
the moving averages of the readings over windows of length $\ell=51$ were
normalized. These normalized averages correspond
to the $Y_{n \ell j}$ scores defined just before
(\ref{Lnlj}). The computed local f.d.r. of the scores at each $j$
determined the conclusion of an aligned signal
there. The scientific purpose is to detect rare
inherited CNV that may occur in a small fraction, perhaps 5\%, of the
population.

Consider, for example, $\ell_T=51$ and
$\pi_N=0.05$. The solution of $T/\ell_T = \exp(N^\zeta-1)$ [see
(\ref{C1})] is $\zeta=0.383$. The solution of $N^{-\beta}=
0.05$ is $\beta=0.568$. Since $1-4 \beta/3 < \zeta
\leq1-\beta$, we are under the second case in (\ref{bNbz}). Based on
(\ref{bNbz}),
the signal-to-noise ratio (for a single
observation in a variant segment) required for successful
detection in a mixture with 5\% variant is then
\[
b_N(\beta,\zeta)/\sqrt{\ell_T} = 0.258.
\]

\begin{figure}

\includegraphics{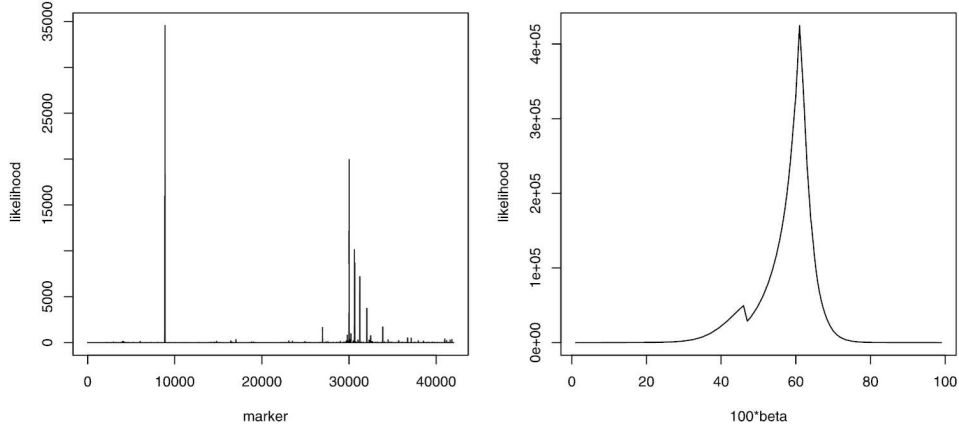}

\caption{(Left) Plot of likelihood against marker position. The
tallest peak is at marker $8852$.
(Right) Plot of likelihood against variant fraction
$\pi_N = N^{-\beta}$ at marker $8852$.}
\label{fig1}
\end{figure}

It is known from earlier studies that the sequence between probes 8800
and 8900
contains two genes that enhance cell death. Copy number losses of these
genes promote unregulated cell
growth, leading to tumor. The display in Figure~\ref{fig1} (left) shows
that the likelihood at marker $j$
(for a signal on the probe interval $j < t \leq j+\ell$),
\[
L_{\ell j} := \int_0^1 \Biggl[ \prod
_{n=1}^N L_{n \ell j}(\beta) \Biggr]
\,d \beta,
\]
is maximized at $j=8852$. The display in Figure~\ref{fig1} (right)
shows that the likelihood at marker
$j=8852$ for variant fraction $\pi_N = N^{-\beta}$,
\[
L_{\ell j}(\beta) := \prod_{n=1}^N
L_{n \ell_T j}(\beta),
\]
is maximized at $\beta=0.61$.
This translates to an estimated 4\% of the population tested having copy
number losses in the probe interval $8852 < t \leq8903$.

%f1 #&#

%s4 #&#
\section{Extensions} \label{extensions}
Cai, Jeng and Jin \cite{CJJ11} and Cai and Wu \cite{CW12} showed that
the HC test is
optimal for heteroscedastic and more general mixtures, respectively.
Arias-Castro and Wang \cite{AW13} analyzed the detection capabilities
of distribution-free tests for
null hypotheses that are not fully specified.
Jeng, Cai and Li \cite{JCL13} showed that the HC test statistic is
optimal for detecting heteroscedastic aligned sparse signals
when assuming that the signal length is very small and that $T$ does
not grow rapidly with $N$.
However, when the aligned signals may range over multiple scales,
the penalty terms introduced in Section~\ref{sec2.2} are critical in ensuring
optimality of the HC test.

Below, we shall show how the detection boundary of Cai et al. \cite
{CJJ11} looks for
general $T/\ell_T$ asymptotics, and show that the adaptive optimality
of the HC and BJ tests extends
to heteroscedastic signals when their penalties, as given in Section~\ref{sec2},
are applied.
This brings home the point that the penalties are not tied down to a
particular model.
Following Jeng et al. \cite{JCL13}, we assume in place of (\ref
{proto}) that
%
%e4.1 #&#
\begin{equation}
\label{X2} X_{nt} = U_{nt} + Z_{nt} \qquad\mbox{where } Z_{nt}
 \mbox{ are i.i.d. N$(0,1)$}.
\end{equation}
Under the null hypothesis $H_0$ of no signals, $U_{nt} \equiv0$ for
all $n$ and $t$.
Under the alternative hypothesis $H_1$ of aligned signals,
there exists an unknown $q>0$ of disjoint intervals $(j_T^{(k)},
j_T^{(k)}+\ell_T^{(k)}]$,
$1 \leq k \leq q$,
such that for the $k$th interval,
%
%e4.2 #&#
\begin{eqnarray}
\label{Unt} U_{nt} & = & \mathrm{ N} \Bigl(\mu_N^{(k)}/
\sqrt{\ell_T^{(k)}}, \tau^{(k)} \Bigr)\qquad
\mbox{if } I_n^{(k)}=1 \mbox{ and } t \in\bigl(j_T^{(k)},
j_T^{(k)}+\ell_T^{(k)}\bigr],
\nonumber
\\[-8pt]
\\[-8pt]
\nonumber
I_n^{(k)} & \sim& \operatorname{ Bernoulli} \bigl(
\pi_N^{(k)} \bigr),
\end{eqnarray}
and $U_{nt}=0$ otherwise, with $\pi_N^{(k)} > 0$, $\mu_N^{(k)} > 0$ and
$\tau^{(k)} \geq0$.
We shall denote $\mu_N^{(1)}$ by $\mu_N$, $\ell_T^{(1)}$ by $\ell_T$
and so forth.
Let $b_N(\beta,\zeta,\tau)$ be such that $b_N(\beta,\zeta,0) =
b_N(\beta
,\zeta)$, and for $\tau> 0$
and $0 \leq\zeta<1-\beta$, let
\begin{eqnarray*}
& & b_N(\beta,\zeta,\tau)/\sqrt{\log N}
\\
&&\qquad =  \cases{
0,  \qquad\mbox{if } \zeta\leq1-2 \beta\mbox{ or } \tau
\geq\displaystyle\frac{\beta}{1-\zeta-\beta},
\vspace*{2pt}\cr
\sqrt{(1-\tau) (2 \beta+\zeta-1)}, \vspace*{2pt}\cr
\hspace*{34pt}\mbox{if } 1-2
\beta< \zeta \leq 1-\displaystyle\frac{4 \beta}{3-\tau},
\vspace*{2pt}\cr
\sqrt{2(1-\zeta)}-\sqrt{2(1+\tau) (1-
\zeta-\beta)},\vspace*{2pt}\cr
\hspace*{34pt}\mbox{if }\displaystyle 1-\min \biggl(2 \beta, \frac{4 \beta}{3-\tau} \biggr) <
\zeta\mbox{ and }\displaystyle \tau< \frac{\beta}{1-\zeta-\beta}.}
\end{eqnarray*}

%th6 #&#
\begin{thmm} \label{thmm6} Assume (\ref{X2}) and (\ref{Unt}).
If for all $1 \leq k \leq q$,
(\ref{C1}) holds and $\mu_N^{(k)} = b_N(\beta^{(k)},\zeta
^{(k)},\tau^{(k)})$, $\pi_N^{(k)} = N^{-\beta^{(k)}-\varepsilon^{(k)}}$
for some $0 < \beta^{(k)} < 1-\zeta^{(k)}$ and $\varepsilon^{(k)} > 0$,
then there is no test that can achieve, at all $j_T^{(k)}$,
$1 \leq k \leq q$,
%
%e4.3 #&#
\begin{equation}
\label{PP0} P(\mbox{Type I error}) + P(\mbox{Type II error}) \rightarrow0.
\end{equation}
Conversely, if for some $1 \leq k \leq q$,
(\ref{C1}) holds and $\mu_N^{(k)} = b_N(\beta^{(k)}, \zeta^{(k)},
\tau^{(k)})$, $\pi_N^{(k)} = N^{-\beta^{(k)}+\varepsilon^{(k)}}$
for some $0 < \beta^{(k)} < 1-\zeta^{(k)}$ and $0 < \varepsilon^{(k)}
\leq
\beta^{(k)}$,
then (\ref{PP0}) can be achieved by the penalized HC and BJ tests.
\end{thmm}

It can be checked that setting $\zeta^{(k)}=0$ will recover for us the boundary
for aligned signals in Jeng et al. \cite{JCL13}. Incidentally, they
assumed that
%
%e4.4 #&#
\begin{equation}
\label{poly2} \log T = o \bigl(N^C \bigr)\qquad \mbox{for all } C > 0,
\end{equation}
which effectively brings us to the case $\zeta^{(k)}=0$.
Corollary~\ref{cor2} below extends the optimality of the HC test in
Jeng et al. \cite{JCL13}
to multiscale signal lengths, by introducing the penalty terms as
described in Section~\ref{sec2}.
In place of (\ref{C1}), let $\zeta_N^{(k)} = \log\log(eT/\ell
_T^{(k)})/\log N$.

%co2 #&#
\begin{cor} \label{cor2} Assume (\ref{X2}) and (\ref{Unt}).
Theorem~\ref{thmm6} holds under (\ref{poly2}) with $\mu
_N^{(k)} = b_N(\beta^{(k)},\zeta_N^{(k)},\tau^{(k)})$ and $0 < \beta
^{(k)} < 1$.
\end{cor}

%s5 #&#
\section{Proofs of Theorems \texorpdfstring{\protect\ref{thmm1}}{1}, 
\texorpdfstring{\protect\ref{thmm2}}{4} and 
\texorpdfstring{\protect\ref{thmm3}}{5}}\label{sec5}

We say that $U_m \stackrel{p}{\sim} V_m$ if $U_m=O_p(V_m)$ and $V_m =
O_p(U_m)$, and that $a_m \gg b_m$ if
$a_m/b_m \rightarrow\infty$. We start with the proof of Theorem~\ref
{thmm1} in Section~\ref{sec5.1}, that detection is asymptotically
impossible below the detection boundary $b_N$, followed by the proofs of
Theorem~\ref{thmm2} (in Section~\ref{sec5.2}) and Theorem~\ref{thmm3} (in
Section~\ref{sec5.3}),
that the average likelihood ratio test is optimal. These proofs are
consolidated in this section as they are
unified by a likelihood ratio approach.
Since the detection problem is easier when $q>1$ compared to $q=1$, we
may assume without loss of generality that $q=1$ under $H_1$
in all the proofs.

%s5.1 #&#
\subsection{Proof of Theorem \texorpdfstring{\protect\ref{thmm1}}{1}}\label{sec5.1}
For $\zeta=0$ the claim of the theorem reduces to theorems proved
by Ingster \cite{Ing98} in the sparse case and by Cai, Jeng and Jin
\cite{CJJ11} in the dense case.

Let $\zeta>0$, and set
$i_T=\lfloor T/\ell_T \rfloor-1$, so $i_T \sim\exp(N^\zeta-1)$, set
$\mu_N = b_N(\beta,\zeta)$ and $\pi_N = N^{-\beta-
\varepsilon}$. Let $Y_{n1} = \ell_T^{-1/2}(X_{n,j_T+1}+\cdots+
X_{n,j_T+\ell_T})$, and let each $Y_{ni}$, $2 \leq i \leq i_T$, be of
the form
$\ell_T^{-1/2}(X_{n,j+1}+\cdots+X_{n,j+\ell_T})$, with all
$(j,j+\ell_T]$
disjoint from each other, and from $(j_T,j_T+\ell_T]$. Let
%
%e5.1 #&#
\begin{equation}
\label{LNi} L_{ni} = 1+\pi_N \bigl[\exp \bigl(
\mu_N Y_{ni}-\mu_N^2/2 \bigr)-1
\bigr],  \qquad L_i = \prod_{n=1}^N
L_{ni}.
\end{equation}
Since $(j,j+\ell_T]$ are disjoint, $L_1, \ldots, L_{i_T}$ are
independent. We take note that
$L_2, \ldots, L_{i_T}$ have
identical distributions which are unchanged when we switch from $H_0$
to $H_1$. In contrast, the distribution of
$L_1$ changes when we switch from $H_0$ to $H_1$. Consider
\[
L = \frac{1}{i_T} L_1 + \frac{1}{i_T} \sum
_{i=2}^{i_T} L_i,
\]
which is the likelihood ratio when $j_T$ is equally likely to
take one of $i_T$ possible values spaced at least $\ell_T$ apart, as explained
above. This
we assume without loss of generality.

If we are able to find $\lambda_N$ such that both
%
%e5.2 #&#
%e5.3 #&#
\begin{eqnarray}
\label{L1} & & L_1=O_p(\lambda_N)\qquad
\mbox{under } H_1\quad \mbox{and}
\\
\label{PaN} & & P \Biggl( a_N + M \lambda_N > \sum
_{i=2}^{i_T} L_i >
a_N \Biggr) \nrightarrow1\qquad \mbox{for all } a_N \in{
\mathbf R} \mbox{ and } M>0
\end{eqnarray}
are satisfied, then $L$ is unable to achieve (\ref{typeI}). If so, then
no test is able to achieve~(\ref{typeI}) because the likelihood ratio
test is the most powerful test.

\textit{Case} 1: $0 < \zeta\leq1-4 \beta/3$, $\mu_N = \sqrt{\log(1+N^{2
\beta
-1+\zeta})}$, $\pi_N = N^{-\beta-\varepsilon}$
with $0 < \varepsilon< \zeta/2$.
Under $H_1$, $Y_{11}, \ldots, Y_{N1}$ are i.i.d. $(1-\pi_N)\mathrm{N}(0,1)+
\pi
_N$N$(\mu_N,1)$. Let
%
%e5.4 #&#
\begin{equation}
\label{E1L} \lambda_N = E_1(L_1) =
\bigl[1+\pi_N^2 \bigl(e^{\mu_N^2}-1 \bigr)
\bigr]^N = \exp \bigl\{ \bigl[1+o(1) \bigr] N^{\zeta-2 \varepsilon} \bigr\}.
\end{equation}
Hence (\ref{L1}) holds.

Let $i \geq2$. Since $Y_{ni} \sim \mathrm{N}(0,1)$,
%
%e5.5 #&#
\begin{equation}
\label{E0V} E_0(L_i)=1,  \qquad E_0
\bigl(L_i^2 \bigr) = \bigl[1+\pi_N^2
\bigl(e^{\mu_N^2}-1 \bigr) \bigr]^N = \lambda_N.
\end{equation}
We check in Appendix \ref{appA} that Lyapunov's condition holds. Hence
%
%e5.6 #&#
\begin{equation}
\label{CLT} \frac{1}{\sqrt{(\lambda_N-1)(i_T-1)}} \sum_{i=2}^{i_T}
(L_i-1) \Rightarrow\mathrm{ N}(0,1).
\end{equation}
Since $\sqrt{(\lambda_N-1)(i_T-1)} = \exp\{[1+o(1)](N^{\zeta-2
\varepsilon
}+N^\zeta)/2 \} \gg\lambda_N$,
(\ref{PaN}) follows from (\ref{CLT}).

\textit{Case} 2: $1-4 \beta/3 < \zeta\leq1-\beta$, $\mu_N = (x-y) \sqrt{2
\log
N}$ where $x=\sqrt{1-\zeta}$, $y=\sqrt{1-\zeta-\beta}$, $\pi_N =
N^{-\beta-\varepsilon}$. Let
\begin{eqnarray*}
\cN^0 & = & \{ n\dvtx I_n=0 \mbox{ or }
Y_{n1} < x \sqrt{2 \log N} \},
\\
\cN^1 & = & \{ n\dvtx I_n=1 \mbox{ and }
Y_{n1} \geq x \sqrt{2 \log N} \},
\end{eqnarray*}
and define, for $h=0,1$,
%
%e5.7 #&#
\begin{equation}
\label{LN1h} L_1^h = \prod
_{n \in\cN^h} L_{n1}\qquad \mbox{where } L_{n1} =
1+\pi_N \bigl[\exp \bigl(\mu_N Y_{n1}-
\mu_N^2/2 \bigr)-1 \bigr].
\end{equation}
Check that
\begin{eqnarray*}
m_0 & := & E_1(L_{n1}|I_n=0) = 1,
\\
m_1 & := & E_1 \bigl[1+(L_{n1}-1) {\mathbf
I}_{\{ Y_{n1} \leq x \sqrt{2 \log N}
\}
}|I_n=1 \bigr]
\\
& =& 1+\pi_N
\bigl[e^{\mu_N^2} \Phi(x \sqrt{2 \log N}-2 \mu_N)-\Phi(x \sqrt{2
\log N}-\mu_N) \bigr].
\end{eqnarray*}
Since $y^2-x^2=-\beta$ and $x < 2(x-y)$ when $1-4 \beta/3 < \zeta
\leq
1-\beta$, it follows that
%
%e5.8 #&#
\begin{eqnarray}
\label{5.6} E_1 \bigl(L_1^0 \bigr) & = &
\bigl[(1-\pi_N)m_0+\pi_N m_1
\bigr]^N\nonumber
\\
& \leq& \bigl[1+\pi_N^2 e^{\mu_N^2} \Phi(x \sqrt{2
\log N}-2 \mu_N) \bigr]^N
\nonumber
\\[-8pt]
\\[-8pt]
\nonumber
& = & \bigl[1+O \bigl(N^{2(y^2-x^2)-2 \varepsilon+2(x-y)^2-(2y-x)^2} \bigr)/ \sqrt{\log N}
\bigr]^N
\\
& = & \exp \bigl[O \bigl(N^{1-x^2-2 \varepsilon} \bigr)/\sqrt{\log N} \bigr] = \exp
\bigl[O \bigl(N^{\zeta-2
\varepsilon} \bigr)/\sqrt{\log N} \bigr].
\nonumber
\end{eqnarray}

Next we apply $\max_{n \in\cN^1} Y_{n1} = O_p(\sqrt{\log N})$ to
show that
%
%e5.9 #&#
\begin{eqnarray}
\label{5.8} \log L_1^1 &=& O_p \bigl(
\bigl(\# \cN^1 \bigr) \log N \bigr)  =  O_p \bigl(N
\pi_N \Phi(-y \sqrt{2 \log N}) \log N \bigr)
\nonumber
\\[-8pt]
\\[-8pt]
\nonumber
& = & O_p \bigl(N^{\zeta-\varepsilon} \sqrt{\log N} \bigr).
\nonumber
\end{eqnarray}
By (\ref{5.6}), (\ref{5.8}) and $L_1 = L_1^0 L_1^1$, (\ref{L1}) holds
for $\lambda_N = \exp(N^{\zeta-\varepsilon}
\log N)$.

Let $i \geq2$. Let $\wtd L_{ni} = 1+\pi_N
[\exp(\mu_N Y_{ni}-\mu_N^2/2)-1] {\mathbf I}_{\{ Y_{ni} \leq x \sqrt{2
\log
N} \}}$ and
$\wtd L_i = \prod_{n=1}^N \wtd L_{ni}$. We check that
%
%e5.10 #&#
\begin{eqnarray}
\label{E0t} E_0(\wtd L_{ni}) & = & 1+\pi_N
\bigl[\Phi(x \sqrt{2 \log N}-\mu_N)-\Phi(x \sqrt{2 \log N}) \bigr]
\nonumber
\\[-8pt]
\\[-8pt]
\nonumber
& = & 1- \bigl[C+o(1) \bigr] N^{-\beta-y^2-\varepsilon}/\sqrt{\log N},
\end{eqnarray}
where $C=(2y \sqrt{\pi})^{-1}$. From this and $1-\beta-y^2=\zeta$, we
conclude that
%
%e5.11 #&#
\begin{equation}
\label{kappaN} \kappa_N := E_0(\wtd L_i) =
\exp \bigl\{ - \bigl[C+o(1) \bigr] N^{\zeta-\varepsilon
}/\sqrt {\log N} \bigr\}.
\end{equation}
Since $E_0(\wtd L_{ni}^2) \geq[E_0(\wtd L_{ni})]^2$ for $n \geq2$, it
follows that
%
%e5.12 #&#
\begin{eqnarray}
\label{vNVar} v_N & := & \operatorname{Var}_0(\wtd
L_i) = \prod_{n=1}^N
E_0 \bigl(\wtd L_{ni}^2 \bigr) -
\kappa_N^2
\nonumber
\\[-8pt]
\\[-8pt]
\nonumber
& \geq& \biggl( \frac{E_0(\wtd L_{1i}^2)}{[E_0(\wtd L_{1i})]^2}-1 \biggr) \kappa_N^2
\sim \operatorname{Var}_0(\wtd L_{1i})
\kappa_N^2,
\end{eqnarray}
and by (\ref{E0t}),
%
%e5.13 #&#
\begin{eqnarray}
\label{Vart} \operatorname{Var}_0(\wtd L_{1i}) & = &
\pi_N^2 \bigl[e^{\mu_N^2} \Phi(x \sqrt{2 \log N}-2
\mu_N)\nonumber
\\
& &\hspace*{16pt}{}  -2 \Phi(x \sqrt{2 \log N}-\mu_N) +\Phi(x \sqrt{2 \log N})
\bigr]
\nonumber
\\[-8pt]
\\[-8pt]
\nonumber
& &{}  - \bigl[E_0(\wtd L_{1i}-1) \bigr]^2
\\
& \sim& C_1 N^{\zeta-1-2 \varepsilon}/\sqrt{\log N},
\nonumber
\end{eqnarray}
where $C_1 = [(2x-4y) \sqrt{\pi}]^{-1}$. We check that Lyapunov's
condition holds and conclude
that
%
%e5.14 #&#
\begin{equation}
\label{Lya2} \frac{1}{\sqrt{v_N (i_T-1)}} \sum_{i=2}^{i_T}
(\wtd L_i-\kappa_N) \Rightarrow\mathrm{ N}(0,1).
\end{equation}
By (\ref{kappaN})--(\ref{Vart}),
$\sqrt{v_N(i_T-1)} \geq\exp\{[1+o(1)] N^\zeta/2 \} \gg\lambda_N$,
and so (\ref{PaN}) holds, but with $L_i$
replaced by $\wtd L_i$. The variability of $L_i$ is larger than that of
$\wtd L_i$, and hence (\ref{PaN}) for
$L_i$ holds as well.

\textit{Case} 3: $\zeta> 1-\beta$, $\mu_N = \sqrt{N^{\zeta-1
+\beta}}$, $\pi_N = N^{-\beta-\varepsilon}$ with $0 < \varepsilon< \zeta
+\beta-1$. Let $\cN^h =
\{ n\dvtx I_n=h \}$, $h=0,1$, and define $L_1^h$ as in (\ref{LN1h}). Since
$\max_{n \in\cN^0} Y_{n1} = O_p(\sqrt{\log N})=o_p
(\mu_N)$, it follows that
%
%e5.15 #&#
\begin{equation}
\label{P1log} P_1 \bigl(L_1^0 \leq1 \bigr)
\rightarrow1.
\end{equation}
We next apply the inequality
\[
\log \bigl(1-\pi_N+\pi_N e^{\mu_N Y_{n1}-
\mu_N^2/2} \bigr) \leq1+
\max \bigl(\mu_N Y_{n1}-\mu_N^2/2 -
\beta\log N,0 \bigr)
\]
to show that
%
%e5.16 #&#
\begin{equation}
\label{logLN11} \log L_1^1 \leq \bigl[1+o_p(1)
\bigr] N \pi_N E \bigl(\mu_N Y_{11}-
\mu_N^2/2|I_1=1 \bigr) \stackrel{p} {\sim}
N^{\zeta-\varepsilon}.
\end{equation}
It follows from (\ref{P1log}), (\ref{logLN11}) and $L_1 = L_1^0 L_1^1$
that (\ref{L1}) holds for
$\lambda_N =  \exp(N^{\zeta-\varepsilon} \log N)$.

Let $i \geq2$ and $\Gamma_i = \{ Y_{ni} \geq\mu_N$ for all $1 \leq n
\leq N^{1-\beta}/2+1 \}$. Then
\[
\log P(\Gamma_i) \sim-N^{1-\beta} \mu_N^2/4
\sim-N^\zeta/2.
\]
Hence $i_T \sim\exp(N^\zeta-1) \gg[P(\Gamma_i)]^{-1}$ and
%
%e5.17 #&#
\begin{equation}
\label{hexi} P \bigl(\#\{ i\dvtx \Gamma_i \mbox{ occurs} \}
=k_N \bigr) \nrightarrow1\qquad \mbox{for any } k_N.
\end{equation}
If $\Gamma_i$ occurs, then
%
%e5.18 #&#
\begin{eqnarray}
\label{Liarray} L_i & \geq& (1-\pi_N)^N
\bigl(1-\pi_N+\pi_N e^{\mu_N^2/2} \bigr)^{N^{1-\beta}/2}\nonumber
\\
& = & \exp \bigl\{ \bigl[1+o(1) \bigr] \bigl[-N^{1-\beta-\varepsilon}+
\bigl(N^{1-\beta}/2 \bigr)N^{\zeta
-1+\beta}/2 \bigr] \bigr\}
\\
& = & \exp \bigl\{ \bigl[1+o(1) \bigr] N^\zeta/4 \bigr\} \gg
\lambda_N.
\nonumber
\end{eqnarray}
Since typically $L_i \ll\lambda_N$, we can conclude (\ref{PaN}) from
(\ref{hexi}) and (\ref{Liarray}).

%s5.2 #&#
\subsection{Proof of Theorem \texorpdfstring{\protect\ref{thmm2}}{4}}\label{sec5.2}

Let (\ref{munt}) and (\ref{C1}) hold with
$\mu_N = b_N(\beta,\zeta)$ and $\pi_N = N^{-\beta+\varepsilon}$. Let
$(j_T^*,j_T^*+\ell_T^*)
\in B_{r,T}$ satisfy (\ref{llr}), and note that by (\ref{C1}) and~(\ref
{BrT}), $r \sim N^\zeta$. Hence for $\zeta>0$ and $N$ large,
%
%e5.19 #&#
\begin{eqnarray}
\label{E1Y} E_1(Y_{n \ell_T^* j_T^*}|I_n=1) & = &
\mu_N \sqrt{\ell_T^*/\ell_T} \geq
\bigl[1-O \bigl(N^{-\zeta/2} \bigr) \bigr] b_N(\beta,\zeta)
\nonumber
\\[-8pt]
\\[-8pt]
\nonumber
& \geq& b_N(\eta,\zeta_{\ell_T^*,NT}), \qquad \beta_1
\leq\eta \leq \beta_2,
\end{eqnarray}
for any $(\beta_1, \beta_2)$ lying in the interior of $(\beta
-\varepsilon
,\beta)$,
and this inequality can also be checked for $\zeta=0$. Let
%
%e5.20 #&#
\begin{equation}
\label{Leta} L(\eta) = \prod_{n=1}^N
L_{n \ell_T^* j_T^*}(\eta).
\end{equation}
Since $r^3 e^{r+1} = O(N^{3 \zeta}
\exp(N^\zeta))$, in view of (\ref{ANT}) and (\ref{E1Y}), Theorem~\ref
{thmm2} follows from
%
%e5.21 #&#
\begin{equation}
\label{check} N^{-3 \zeta} \exp \bigl(-N^\zeta \bigr) \int
_{\beta_1}^{\beta_2} L(\eta) \,d \eta \stackrel{p} {
\rightarrow} \infty.
\end{equation}
To show (\ref{check}), it suffices to check that
%
%e5.22 #&#
\begin{equation}
\label{check2} N^{-3 \zeta} \exp \bigl(-N^\zeta \bigr) L(\eta)
\stackrel{p} {\rightarrow} \infty,
\end{equation}
when $E_1(Y_{n \ell_T^* j_T^*}|I_n=1) \geq b_N(\eta,\zeta_{\ell
_T^*,NT})$, and $\pi_N = N^{-\beta+\varepsilon}$
with $\beta-\varepsilon< \eta$.
To see why (\ref{check2}) leads to (\ref{check}), define $C_{N \zeta
}=N^{3 \zeta} \exp(N^\zeta)$, let
$M > 0$ and let
\[
\chi_N = \mbox{Leb. meas.} \bigl\{ \eta\in[\beta_1,
\beta_2]\dvtx L(\eta) < M C_{N \zeta} \bigr\}.
\]
By (\ref{check2}), $E_1 \chi_N \rightarrow0$ and hence $P_1 \{ \chi_N
< (\beta_2-\beta_1)/2
\} \rightarrow1$. If $\chi_N < (\beta_2-\beta_1)/2$, then $C_{N
\zeta
}^{-1} \int_{\beta_1}^{\beta_2}
L(\eta) \,d \eta\geq M (\beta_2-\beta_1)/2$. Since $M > 0$ can be chosen
arbitrarily large, (\ref{check}) holds.

To cross-reference the results in the proof of Theorem~\ref{thmm1} more
easily, we relabel $j_T^*$ and $\ell_T^*$
in (\ref{Leta}) by $j_T$ and $\ell_T$, respectively, and rephrase
(\ref
{check2}) as
%
%e5.23 #&#
\begin{equation}
\label{N3zeta} N^{-3 \zeta} \exp \bigl(-N^\zeta \bigr) L(\beta)
\stackrel{p} {\rightarrow} \infty,
\end{equation}
when $\mu_N = b_N(\beta,\zeta)$ and $\pi_N = N^{-\beta+\varepsilon}$.

\textit{Case} 1: $0 \leq\zeta\leq1-4 \beta/3$,\vspace*{1pt} $\mu_N = \sqrt{\log(1+N^{2
\beta-1+\zeta})}$, $\pi_N = N^{-\beta+\varepsilon}$
with $0 < \varepsilon< (1-\zeta)/2$.
Under $H_1$, $Y_{n1} = \mu_N {\mathbf I}_{\{ I_n=1 \}}+Z_{n1}$, where
$Z_{11}, \ldots, Z_{N1}$ are i.i.d. $\mathrm{N}(0,1)$.
Let $L_{n1} = 1+N^{-\beta} [\exp(\mu_N Z_{n1}-\mu_N^2/2)-1]$, and
%
%e5.24 #&#
%e5.25 #&#
\begin{eqnarray}
\label{L10} L_1^0 & = & \prod
_{n=1}^N L_{n1},
\\
\label{L11} L_1^1 & =& \prod
_{n: I_n=1} \biggl( \frac{1+N^{-\beta} [\exp(\mu_N
Y_{n1}-\mu_N^2/2)-1]}{ 1+N^{-\beta}
[\exp(\mu_N Z_{n1}-\mu_N^2/2)-1]} \biggr).
\end{eqnarray}

Since for $v \geq0$,
%
%e5.26 #&#
\begin{equation}
\label{vincrease} f(v) := \frac{1+N^{-\beta} (ve^{\mu_N^2}-1)}{1+N^{-\beta} (v-1)}
\qquad\mbox{is increasing and } f(v) \geq1,
\end{equation}
it follows that
%
%e5.27 #&#
\begin{eqnarray}
\label{logL11}  \log L_1^1& \geq&\# \{ n\dvtx
I_n=1, Y_{n1} \geq2 \mu_N \} \log \biggl(
\frac{1+N^{-\beta}
(e^{3 \mu_N^2/2}-1)}{1+N^{-\beta} (e^{\mu_N^2/2}-1)} \biggr)
\nonumber
\\[-8pt]
\\[-8pt]
\nonumber
& \stackrel{p} {\sim} & N \pi_N \Phi(-\mu_N) \times
\cases{ %
N^{\beta-1+\zeta}, &\quad $\mbox{if } 0 \leq\zeta<
1-2 \beta,$
\vspace*{2pt}\cr
N^{-\beta} \sqrt{2}, &\quad $\mbox{if } \zeta =1-2 \beta,$
\vspace*{2pt}\cr
N^{-\beta} e^{3 \mu_N^2/2}, &\quad $\mbox{if } 1-2 \beta< \zeta< 1-4 \beta /3,$
\vspace*{2pt}\cr
\log2, &\quad $\mbox{if } \zeta=1-4 \beta/3.$}
\end{eqnarray}
In the above, we apply the relation
$\log[1+N^{-\beta}(e^{\kappa\mu_N^2/2}-1)] \sim N^{-\beta}
(e^{\kappa
\mu_N^2/2}-1)$
for $\kappa=1,3$, with the exception
\[
\log \bigl[1+N^{-\beta} \bigl(e^{3 \mu_N^2/2}-1 \bigr) \bigr] \sim\log2\qquad
\mbox{when } \zeta =1-4 \beta/3.
\]

Since $\Phi(-\mu_N) \sim\frac{1}{2}$ when $0 \leq\zeta< 1-2 \beta$,
$\Phi(-\mu_N) = \Phi(-\sqrt{\log2})$ when
$\zeta=1-2 \beta$ and $\Phi(-\mu_N) \sim e^{-\mu_N^2/2}/(\mu_N
\sqrt{2
\pi})$ when $1-2 \beta< \zeta\leq
1-4 \beta/3$, it follows from checking each of the cases in (\ref
{logL11}) that
%
%e5.28 #&#
\begin{equation}
\label{logL11N} \frac{\log L_1^1}{N^{\zeta+\varepsilon}(\log N)^{-1}} \stackrel {p} {\rightarrow} \infty.
\end{equation}

We shall next obtain lower bounds of $\log L_1^0$. We apply Taylor's
expansion $\log(1+u)=u-[\frac{1}{2}+o(1)]
u^2$ to show that
%
%e5.29 #&#
\begin{equation}
\label{ElogLn1} E_1(\log L_{n1}) \sim-N^{-2 \beta}
\bigl(e^{\mu_N^2}-1 \bigr)/2 = -N^{\zeta-1}/2,
\end{equation}
and $\log(1+u) \sim u$ to show that
%
%e5.30 #&#
\begin{equation}
\label{ElogLn12} E_1(\log L_{n1})^2 \sim
N^{-2 \beta} \bigl(e^{\mu_N^2}-1 \bigr) = N^{\zeta-1}.
\end{equation}
It follows from (\ref{ElogLn1}) and (\ref{ElogLn12}) that $\log L_1^0 =
\sum_{n=1}^N \log L_{n1}
\stackrel{p}{\sim} -N^\zeta/2$.
Since $L(\beta) = L_1^0 L_1^1$, we can conclude (\ref{N3zeta}) from
(\ref{logL11N}).

\textit{Case} 2: $1-4 \beta/3 < \zeta\leq1-\beta$, $\mu_N= (x-y) \sqrt{2
\log
N}$ where $x=\sqrt{1-\zeta}$, $y=\sqrt{1-
\zeta-\beta}$, $\pi_N = N^{-\beta+\varepsilon}$. Define
\begin{eqnarray*}
\wtd L_{n1} & = & 1+N^{-\beta} \bigl[\exp \bigl(
\mu_N Z_{n1}-\mu_N^2/2 \bigr)-1
\bigr] {\mathbf I}_{\{ Z_{n1} \leq
x \sqrt{2 \log N} \}},
\\
\wtd L_1^0 & = &
\prod_{n=1}^N \wtd L_{n1},
\\
\wtd L_1^1 & =& \prod_{n: I_n=1}
\biggl( \frac{1+N^{-\beta} [\exp
(\mu_N
Y_{n1}-\mu_N^2/2)-1]
{\mathbf I}_{\{ Z_{n1} \leq
x \sqrt{2 \log N} \}}}{ 1+N^{-\beta} [\exp(\mu_N Z_{n1}-\mu
_N^2/2)-1]{\mathbf I}_{\{ Z_{n1} \leq
x \sqrt{2 \log N} \}}
} \biggr).
\end{eqnarray*}
By (\ref{vincrease}),
%
%e5.31 #&#
\begin{eqnarray}
\label{loghex}  \log\wtd L_1^1& \geq&\# \{ n\dvtx
I_n=1, x \sqrt{2 \log N} \leq Y_{n1} \leq2x \sqrt{2 \log
N} \}
\nonumber\\
& &{}  \times\log \biggl( \frac{1+N^{-\beta}[\exp(\mu_N x \sqrt{2
\log N} -
\mu_N^2/2)-1]}{1+N^{-\beta} [\exp(\mu_N x \sqrt{2 \log N}-3 \mu
_N^2/2)-1]} \biggr)
\\
& \sim& N \pi_N \Phi(-y \sqrt{2 \log N}) \log2 \sim C (\log2)
N^{\zeta+\varepsilon}/\sqrt{\log N},
\nonumber
\end{eqnarray}
where $C=(2y \sqrt{\pi})^{-1}$. Recall that $C_1=[(2x-4y) \sqrt{\pi}]^{-1}$.

Apply Taylor's expansion $\log(1+u)=u-[\frac{1}{2}+o(1)] u^2$ to show that
%
%e5.32 #&#
\begin{eqnarray}
\label{Elogt} & & E_1(\log\wtd L_{n1})\nonumber
\\
&&\qquad =  N^{-\beta} \bigl[\Phi(y \sqrt{2 \log N})-\Phi(x \sqrt{2 \log N})
\bigr]
\nonumber
\\
& &\qquad\quad{}  - \bigl[1/2+o(1) \bigr] N^{-2 \beta} \bigl[e^{\mu_N^2} \Phi
\bigl((2y-x) \sqrt{2 \log N} \bigr)
\nonumber
\\[-8pt]
\\[-8pt]
\nonumber
& &\hspace*{98pt}\qquad\quad{}  -2 \Phi(y \sqrt{2 \log N})+\Phi(x \sqrt{2 \log N}) \bigr]
\\
&&\qquad =  - \bigl[1+o(1) \bigr] CN^{-\beta-y^2}/\sqrt{\log N}- \bigl[1/2+o(1)
\bigr] C_1 N^{-x^2}/\sqrt{\log N}
\nonumber
\\
&&\qquad =  - \bigl[1+o(1) \bigr] (C+C_1/2) N^{\zeta-1}/\sqrt{\log
N},
\nonumber
\end{eqnarray}
and $\log(1+u) \sim u$ to show that
%
%e5.33 #&#
\begin{equation}
\label{Elogt2} E_1(\log\wtd L_{n1})^2 \sim
C_1 N^{\zeta-1}/\sqrt{\log N}.
\end{equation}
It follows from (\ref{Elogt}) and (\ref{Elogt2}) that
\[
\log\wtd L_1^0 = \sum_{n=1}^N
\log\wtd L_{n1} \stackrel{p} {\sim} -(C+C_1/2)
N^\zeta/\sqrt{\log N}
\]
if $\zeta>0$ and $|\log\wtd L_1^0| =O_p(1)$ if $\zeta=0$.
Since $L(\beta) \geq\wtd L_1^0 \wtd L_1^1$, we can conclude~(\ref
{N3zeta}) from (\ref{loghex}).

\textit{Case} 3: $\zeta> 1-\beta$, $\mu_N = \sqrt{N^{\zeta-1+\beta}}$,
$\pi_N =
N^{-\beta+\varepsilon}$.
The inequality
\[
L(\beta) \geq \bigl(1-N^{-\beta} \bigr)^N \prod
_{n:I_n=1} \bigl\{ 1+N^{-\beta} \bigl[\exp \bigl(\mu
_N Y_{n1}-\mu_N^2/2 \bigr)-1
\bigr] \bigr\}
\]
leads to
\begin{eqnarray*}
\log L(\beta) & \geq& -2N^{1-\beta}+ \bigl[1+o_p(1) \bigr] N
\pi_N \bigl[\mu_N E_1(Y_{n1})-
\mu_N^2/2-\beta\log N \bigr]
\\
& = & -2N^{1-\beta}+
\bigl[1+o_p(1) \bigr] N^{\zeta+\varepsilon}/2 \stackrel {p} {\sim}
N^{\zeta+\varepsilon}/2,
\end{eqnarray*}
and from this, we can conclude (\ref{N3zeta}).

%s5.3 #&#
\subsection{Proof of Theorem \texorpdfstring{\protect\ref{thmm3}}{5}}\label{sec5.3}

Let
\[
E_1(Y_{\ell_T j_T})=b_T(\ell_T)+c_T\qquad
\mbox{with } c_T \rightarrow \infty,
\]
and let $(j_T^*,j_T^*+\ell_T^*) \in B_{r,T}$ be such that (\ref{llr})
holds. Hence
%
%e5.34 #&#
\begin{equation}\qquad
\label{E1Ystar} E_1(Y_{\ell_T^* j_T^*}) = \bigl[b_T(
\ell_T)+c_T \bigr] \sqrt{
\ell_T^*/ \ell_T} = \bigl[1-O \bigl(r^{-1/2}
\bigr) \bigr] \bigl[b_T(\ell_T)+c_T
\bigr].
\end{equation}
Since $b_T(\ell_T^*) = \sqrt{2 \log(eT/\ell_T^*)} = \sqrt
{b_T^2(\ell
_T)+O(1)} = b_T(\ell_T)+O(1)$
and\break  $r^{-1/2} b_T(\ell_T) = O(b_T(\ell_T)/\sqrt{\log(T/\ell_T)}) =
O(1)$, it follows from
(\ref{E1Ystar}) that
%
%e5.35 #&#
\begin{equation}
\label{cTprime} E_1(Y_{\ell_T^* j_T^*}) = b_T \bigl(
\ell_T^* \bigr)+c_T' \qquad\mbox{with }
c_T' \rightarrow\infty.
\end{equation}
We check that under (\ref{cTprime}),
\[
\frac{\exp[b_T(\ell_T^*) Y_{\ell_T^* j_T^*}-b_T^2(\ell
_T^*)/2]}{[\log
(T/\ell_T^*)]^3 (T/\ell_T^*)} \stackrel{p} {\rightarrow} \infty.
\]
Since $r^3 e^{r+1} = O([\log(T/\ell_T)]^3(T/\ell_T))$, it follows from
(\ref{AT}) that $A_T \stackrel{p}{
\rightarrow} \infty$.

%s6 #&#
\section{Proofs of Theorems 
\texorpdfstring{\protect\ref{thmm2a}}{2} and 
\texorpdfstring{\protect\ref{thmm2b}}{3}}\label{sec6}

We shall prove Theorem~\ref{thmm2b} in Section~\ref{sec6.1}, that the penalized
Berk--Jones test is optimal, and Theorem~\ref{thmm2a} in
Section~\ref{sec6.2}, that the penalized higher criticism test
is optimal as well.

%s6.1 #&#
\subsection{Proof of Theorem \texorpdfstring{\protect\ref{thmm2b}}{}}\label{sec6.1}

In Lemma~\ref{lem1} below, we show that the
Type I error probability of the
penalized Berk--Jones test statistic goes
to zero for the threshold $h_N := 2 \log N$. We do this in more
generality than is required for proving
Theorem~\ref{thmm2b}. For each $\ell$ and $j$, we assume only that $p_{1
\ell j}, \ldots, p_{N \ell j}$ are i.i.d. $\operatorname{Uniform}(0,1)$
random variables under $H_0$. Hence we allow for $X_{nt}$ to be
non-Gaussian, and for $X_{n1}, \ldots, X_{nT}$ to be dependent
random variables under $H_0$.

%le1 #&#
\begin{lem} \label{lem1} Assume that for each $\ell$ and $j$, $p_{(1)
\ell j} \leq
\cdots\leq p_{(N) \ell j}$ in (\ref{BJNlj}) and (\ref
{PBJ}) are the ordered values of i.i.d.  $\operatorname{Uniform}(0,1)$ random
variables. Then
%
%e6.1 #&#
\begin{equation}
\label{l1} P_0 \{ \mathrm{ PBJ}_{NT} \geq
h_N \} \rightarrow0 \qquad\mbox{as } N \rightarrow \infty.
\end{equation}
\end{lem}

\begin{pf}Let $a_{N \ell} = h_N + s_{\ell T} \log s_{\ell T}$.
For each $1 \leq\ell\leq T$ and $1 \leq n \leq N$, let $\rho_{\ell}$
be such that $K(\frac{n}{N},\rho_{\ell})=a_{N \ell}/N$. Let
$\bar S_N(t)\ [=\bar
S_{N \ell j}(t)]= N^{-1}\times\break  \sum_{n=1}^N {\mathbf I}_{ \{ Y_{n \ell j} \geq
z(t) \}}$, where $z(t)$ denotes the upper $t$-quantile
of the standard normal. By the Chernoff--Hoeffding inequality,
\begin{eqnarray}
\label{PSN} & & P_0 \bigl\{ K(n/N, p_{(n)}) \geq
a_{N \ell}/N \bigr\}\nonumber
\nonumber\\
\eqntext{ \bigl(=  P_0 \bigl\{ \bar S_N(
\rho_{\ell}) \geq n/N \bigr\} \bigr) \leq e^{-a_{N \ell}} =
N^{-2} e^{-s_{\ell T} \log s_{\ell T}},}
\end{eqnarray}
and hence by Bonferroni's inequality,
%
%e6.2 #&#
\begin{equation}
\label{BJcs} P_0 \{ \mathrm{ BJ}_{N \ell j} -
s_{\ell T} \log s_{\ell T} \geq h_N \} \leq
N^{-1} e^{-s_{\ell T} \log s_{\ell T}}.
\end{equation}
By (\ref{BrT}), $\# B_{r,T} \leq re^{r+1}$.
Since $\ell\leq T/e^{r-1}$ for $(j,j+\ell) \in B_{r,T}$, so
$s_{\ell T} = \log(eT/\ell) \geq r$, and by (\ref{PBJ}) and
(\ref{BJcs}),
%
%e6.3 #&#
\begin{equation}
\label{6.36a} P_0 \{ \mathrm{PBJ}_{NT} \geq
h_N \} \leq N^{-1} \sum_{r=1}^{\infty}
r e^{r+1-r \log r},
\end{equation}
and (\ref{l1}) holds.
\end{pf} %$\wbox$

\begin{pf*}{Proof of Theorem \ref{thmm2b}} Let
$(j^*,j^*+\ell^*)\ [=(j_T^*,j_T^*+\ell_T^*)]
\in B_{r,T}$ be such that $j_T \leq j^* < j^*+\ell^*
\leq j_T+ \ell_T$ and $1-\ell^*/\ell_T =
O(r^{-1/2})$; see (\ref{llr}). Since $a_{N \ell^*} = O(N^\zeta\log N)$,
in view of Lemma~\ref{lem1},
it remains for us to show that if $\mu_N = b_N(\beta,\zeta)$ and
$\pi_N =
N^{-\beta+\varepsilon}$ for $\varepsilon> 0$, then we can find $t_N$ such
that in each case below,
%
%e6.4 #&#
\begin{equation}
\label{PK} \frac{K(\bar S_N,t_N)}{N^{\zeta-1} \log N}
\rightarrow\infty\qquad \mbox{where } \bar
S_N = \bar S_{N \ell^* j^*}(t_N).
\end{equation}

\textit{Case} 1(a): $0 \leq\zeta\leq1-2 \beta$,
$b_N(\beta,\zeta) = \sqrt{\log(1+N^{2\beta-1+\zeta})}$.
Let $t_N =\break  \Phi(-2 b_N(\beta,\zeta))$. Except when
$\zeta=1-2 \beta$, we have $b_N(\beta,\zeta) \rightarrow0$ and
%
%e6.5 #&#
\begin{equation}
\label{E1S} E_1 \bar S_N - t_N \sim(2
\pi)^{-1/2} \pi_N b_N(\beta,\zeta) \sim(2
\pi)^{-1/2} N^{(\zeta-1)/2+\varepsilon}.
\end{equation}
By Taylor's expansion, $K(t,x) \sim2(t-x)^2$ when $t \rightarrow\frac
{1}{2}$ and $x \rightarrow\frac{1}{2}$. Moreover,
the standard error of $\bar S_N\ [\sim
(4N)^{-1/2}]$ is small relative to (\ref{E1S}). Hence (\ref{PK}) holds
because
%
%e6.6 #&#
\begin{equation}
\label{1aK} K ( \bar S_N, t_N ) \stackrel{p} {\sim} 2
( \bar S_N - t_N)^2 \stackrel {p} {\sim}
\pi^{-1} N^{\zeta-1+2 \varepsilon}.
\end{equation}

When $\zeta=1-2 \beta$, $b_N(\beta,\zeta) = \sqrt{\log2}$ and
\[
E_1 \bar S_N - t_N \sim\wtd C
\pi_N\qquad \mbox{where } \wtd C = \Phi(-2 \sqrt{\log2})- \Phi(-\sqrt{
\log2}).
\]
This leads to (\ref{1aK})
with $\pi^{-1}$ replaced by $2 \wtd C^2$, and then to
(\ref{PK}).

\textit{Case} 1(b): $1-2 \beta< \zeta< 1-4 \beta/3$,
$b_N(\beta,\zeta) = \sqrt{\log(1+N^{2 \beta-1+\zeta})}$\break 
($\sim C \sqrt{\log N}$, where $C = \sqrt{2 \beta-1+\zeta}$).
Let $t_N = \Phi(-2 b_N(\beta,\zeta))$. Then
%
%e6.7 #&#
%e6.8 #&#
%e6.9 #&#
\begin{eqnarray}
\label{1btN} t_N & \sim& (C \sqrt{8 \pi\log N})^{-1}
N^{-4 \beta-2(\zeta-1)},
\\
\label{1bdiff} E_1 \bar S_N - t_N & \sim &
\pi_N \Phi \bigl(-b_N(\beta,\zeta) \bigr)
\nonumber
\\[-8pt]
\\[-8pt]
\nonumber
& \sim& (C \sqrt{2 \pi\log N})^{-1} N^{-2 \beta-(\zeta
-1)/2+\varepsilon
},
\\
\label{1bV} \operatorname{Var}_1 \bar S_N & \sim&
N^{-1} \bigl[t_N+\pi_N \Phi
\bigl(-b_N( \beta,\zeta) \bigr) \bigr].
\end{eqnarray}

We claim that a consequence of (\ref{1btN})--(\ref{1bV}) is that
%
%e6.10 #&#
\begin{equation}
\label{1bdiff2} \sqrt{t_N N^{\zeta-1} \log N} + \sqrt{
\operatorname{Var}_1 \bar S_N} = o(E_1 \bar
S_N -t_N).
\end{equation}
By (\ref{1bdiff2}), $\sqrt{t_N N^{\zeta-1} \log N}
= o_p(|\bar S_N-t_N|)$ and hence
%
%e6.11 #&#
\begin{equation}
\label{1bgg} \frac{(\bar S_N-t_N)^2/(2t_N)}{N^{\zeta-1} \log N} \stackrel{p} {\rightarrow} \infty.
\end{equation}

By (\ref{1btN}), the solution in $y$ of $y^2/(2t_N) = a_{N \ell^*}/N$
satisfies
\[
y \sim \biggl( \frac{\zeta^*}{C} \sqrt{\frac{\log N}{2 \pi}}
\biggr)^{1/2} N^{-2 \beta-(\zeta-1)/2} \qquad\bigl[=o(t_N) \mbox{ because }
\zeta< 1-4 \beta/3 \bigr],
\]
where $\zeta^* = \zeta$ if $\zeta\neq0$ and $\zeta^*=2$ if $\zeta=0$.
Hence by (\ref{1bgg}) and $K(x,t) \sim\frac{(x-t)^2}{2t}$, as $t
\rightarrow
0$ and $\frac{x}{t} \rightarrow1$, (\ref{PK}) holds.

It remains for us to show (\ref{1bdiff2}). By (\ref{1btN}), the
exponent of $N$ in
$\sqrt{t_N N^{\zeta-1}}$
is $-2 \beta-(\zeta-1)/2$, which is smaller than the exponent in $N$ of
$E_1 \bar S_N-t_N$;
see~(\ref{1bdiff}). Therefore,
%
%e6.12 #&#
\begin{equation}
\label{6.46} \sqrt{t_N N^{\zeta-1} \log N} =
o(E_1 \bar S_N-t_N).
\end{equation}
The leading exponent of $N$ in Var$_1 \bar S_N$ is
\[
\max \bigl(-4 \beta-2 \zeta+1,-2 \beta-(\zeta+1)/2+\varepsilon \bigr)\qquad \bigl[< -4
\beta-(\zeta-1)+ 2 \varepsilon \bigr],
\]
and therefore by (\ref{1bdiff}), Var$_1 S_N = o((E_1 \bar
S_N-t_N)^2)$. This,
together with (\ref{6.46}), implies (\ref{1bdiff2}).

\textit{Case} 2: $1-\frac{4}{3} \beta< \zeta\leq1-\beta$,
$b_N(\beta,\zeta) = (x-y) \sqrt{2 \log N}$, where $x=\sqrt{1-\zeta}$,
$y=\sqrt{1-\beta-\zeta}$.
Let $t_N = \Phi(-x \sqrt{2 \log N})\ [\sim(2x \sqrt{\pi\log N})^{-1}
N^{\zeta-1}]$. Then
\[
E_1 \bar S_N \sim(1-\pi_N) t_N +
\pi_N \Phi(-y \sqrt{2 \log N}) \sim (2y \sqrt{\pi\log
N})^{-1} N^{\zeta-1+\varepsilon},
\]
which is large relative to $t_N$, and
\begin{eqnarray*}
\operatorname{Var}_1 \bar S_N & \sim& N^{-1}
\bigl[\Phi(-x \sqrt{2 \log N})+\pi_N \Phi(-y \sqrt{2 \log N}) \bigr]
\\
& \sim& (2y \sqrt{\pi\log N})^{-1} N^{\zeta-2+\varepsilon} = o
\bigl((E_1 \bar S_N)^2 \bigr).
\end{eqnarray*}
Therefore $\bar S_N \stackrel{p}{\sim}
(2y \sqrt{\log N})^{-1} N^{\zeta-1+\varepsilon}$, and by $K(x,t) \sim x
\log\frac{x}{t}$, as
$x \rightarrow0$ and $\frac{x}{t} \rightarrow\infty$,
\[
K(\bar S_N,t_N) \stackrel{p} {\sim} \bar S_N
\log(\bar S_N/t_N) \stackrel{p} {\sim} C' (
\log N)^{1/2} N^{\zeta-1+\varepsilon},
\]
for some $C' > 0$, and (\ref{PK}) therefore holds.

\textit{Case} 3: $\zeta> 1-\beta$, $b_N(\beta,\zeta) = \sqrt{N^{\beta
+\zeta
-1}}$. Let $t_N = \Phi(-b_N(\beta,\zeta)/2)$. Then
\[
t_N \stackrel{p} {\sim} \bigl(\pi N^{\beta+\zeta-1}/2
\bigr)^{-1/2} \exp \bigl(-N^{\beta
+\zeta-1}/8 \bigr) \quad\mbox{and}\quad \bar
S_N \stackrel{p} {\sim} \pi_N = N^{-\beta+\varepsilon}.
\]
Therefore $\bar S_N/t_N \stackrel{p}{\rightarrow} \infty$ and by
$K(x,t) \sim x \log\frac{x}{t}$, as $x \rightarrow0$ and
$\frac{x}{t} \rightarrow\infty$,
\[
K ( \bar S_N, t_N ) \stackrel{p} {\sim} \bar
S_N \log(\bar S_N/t_N) \stackrel{p} {\sim}
N^{\zeta-1+\varepsilon}/8,
\]
and (\ref{PK}) therefore holds.
\end{pf*}

%s6.2 #&#
\subsection{Proof of Theorem \texorpdfstring{\protect\ref{thmm2a}}{}}\label{sec6.2}

In Lemma~\ref{lem2} below, we show that the Type I error probability of
the penalized higher criticism test statistic
goes to zero for the threshold $h_N = 2 \log N$. Again as in Lemma~\ref
{lem1}, we do this more generally than is required for
proving Theorem~\ref{thmm2a}.

%le2 #&#
\begin{lem} \label{lem2} Assume that for each $\ell$ and $j$, $p_{(1)
\ell j} \leq
\cdots\leq p_{(N) \ell j}$ in (\ref{HCNlj}) and (\ref
{PHC}) are the ordered values of i.i.d. $\operatorname{Uniform}(0,1)$ random
variables. Then
\[
P_0 \{ \mathrm{ PHC}_{NT} \geq h_N \}
\rightarrow0 \qquad\mbox{as } N \rightarrow\infty.
\]
\end{lem}

\begin{pf} We first modify (\ref{PSN})--(\ref{6.36a}),
in the proof of Lemma~\ref{lem1}, step-by-step to show that for
$c_{N \ell} := 2 \log N + s_{\ell T} + 3 \log s_{\ell T}$,
\[
P_0 \bigl\{ \mathrm{BJ}_{N \ell j} \geq c_{N \ell} \mbox{
for some } (j,j+\ell ) \in B_T \bigr\} \rightarrow0.
\]
We then combine this with (\ref{KQ2}) in Appendix \ref{appB} to show that
%
%e6.13 #&#
\begin{eqnarray}
\label{Pa}&& P_0 \bigl\{ \mathrm{HC}_{N \ell j} \geq(6
c_{N \ell})^{1/2} \bigl(4+s_{\ell
T}^{-1} \log N
\bigr)^{1/4} \mbox{ for some } (j,j+\ell) \in B_T \bigr\}
\nonumber
\\[-8pt]
\\[-8pt]
\nonumber
&&\qquad
\rightarrow0.
\end{eqnarray}
Therefore it suffices to show that
%
%e6.14 #&#
\begin{equation}\qquad
\label{6.47} (6 c_{N \ell})^{1/2} \bigl(4+s_{\ell T}^{-1}
\log N \bigr)^{1/4} \leq h_N + (s_{\ell T} \log
s_{\ell T})^{1/2} \qquad\mbox{for } N \mbox{ large},
\end{equation}
uniformly over $1 \leq\ell\leq T$. The left-hand side of (\ref{6.47})
is $o(h_N)$ uniformly over $s_{\ell T} \leq\log N$, and
$o((s_{\ell T} \log s_{\ell T})^{1/2})$
uniformly over $s_{\ell T} > \log N$,
so (\ref{6.47}) indeed holds.
\end{pf} %$\wbox$

\begin{pf*}{Proof of Theorem \ref{thmm2a}}
We apply the proofs for
cases 1 and 2 in Theorem~\ref{thmm2b} to show that
there exist $(j^*,j^*+\ell^*) \in
B_T$ and $t_N^* \in[\frac{s_{\ell^* T}}{N}, \frac{1}{2}]$
[$t_N^* = t_N$ for case 1 and $t_N^* =
\frac{s_{\ell^* T}}{N}\ (\sim N^{\zeta-1})$ for case 2] such that
\[
P_1 \bigl\{ K \bigl( \bar S_{N \ell^* j^*} \bigl(t_N^*
\bigr),t_N^* \bigr) \geq N^{\zeta-1+\varepsilon/2} \bigr\} \rightarrow1.
\]
Hence by $\frac{x-t}{\sqrt{t(1-t)}} \geq\sqrt{2 K(x,t)}$ for
$x \geq t$,
\[
P_1 \bigl\{ \mathrm{PHC}_{NT} \geq\sqrt{2N^{\zeta+\varepsilon/2}}-
\sqrt{s_{\ell^* T} \log s_{\ell^* T}} \bigr\} \rightarrow1.
\]
Since $h_N+\sqrt{s_{\ell^* T} \log s_{\ell^* T}} = o(\sqrt{N^{\zeta
+\varepsilon/2}})$, the Type II error
probability indeed
goes to zero.
\end{pf*} %$\wbox$

%s7 #&#
\section{Proofs of Theorem \texorpdfstring{\protect\ref{thmm6}}{6} and 
Corollary
\texorpdfstring{\protect\ref{cor2}}{2}}\label{sec7}

We prove Theorem~\ref{thmm6} here and in Sections \ref{sec7.1} and \ref{sec7.2}.
In Section~\ref{sec7.3}, we prove Corollary~\ref{cor2}.
Let $\mu_N = b_N(\beta,\zeta,\tau)$, $Y_{n1} = \ell_T^{-1/2} \sum_{\ell
=1}^{\ell_T} X_{n,j_T+ \ell}$,
and for $2 \leq i \leq i_T\ (=\lfloor T/\ell_T \rfloor-1)$, let $Y_{ni}
= \ell_T^{-1/2} \sum_{\ell=1}^{\ell_T} X_{n,j+\ell}$,
with all $(j,j+\ell_T]$ disjoint from each other, and from
$(j_T,j_T+\ell_T]$.
Let $L_i = \prod_{n=1}^N L_{ni}$, where
%
%e7.1 #&#
\begin{equation}
\label{Lni6} L_{ni} = 1+\pi_N \biggl\{
\frac{1}{\sqrt{1+\tau}} \exp \biggl[- \frac
{(Y_{ni}-\mu_N)^2}{2(1+\tau)}+ \frac{Y_{ni}^2}{2} \biggr]-1
\biggr\}
\end{equation}
is the likelihood ratio of $Y_{ni} \sim(1-\pi_N) \mathrm{ N}(0,1)+\pi_N
\mathrm{ N}(\mu_N,1+\tau)$ and
$Y_{ni} \sim\mathrm{ N}(0,1)$.
Below, we go over the relevant cases to show that there exists
$\lambda_N$ satisfying (\ref{L1}) and (\ref{PaN}) when $\pi_N =
N^{-\beta-\varepsilon}$.
This implies that there is no test able to achieve (\ref{PP0}).
We shall only consider $\tau> 0$ as the case $\tau=0$ has been
covered in
Theorem~\ref{thmm1}.

\textit{Case} 1: $1-2 \beta< \zeta\leq1-\frac{4 \beta}{3-\tau}$
($\Rightarrow
\tau<1$), $\mu_N = \sqrt{(1-\tau)(2 \beta+\zeta-1) \log N}$.
Let $C = (1-\tau^2)^{-1/2}$. By (\ref{Lni6}),
\[
E_1(L_{n1}) = 1+\pi_N^2
\bigl(Ce^{\mu_N^2/(1-\tau)}-1 \bigr) =1+ \bigl[C+o(1) \bigr] N^{\zeta
-1-2 \varepsilon},
\]
and therefore (\ref{L1}) holds with $\lambda_N = E_1(L_1)\ (= \exp\{
[C+o(1)] N^{\zeta-2 \varepsilon} \})$.
For $i \geq2$, $Y_{ni} \sim \mathrm{N}(0,1)$, $E_0(L_i)=1$ and
$E_0(L_i^2)=E_1(L_1)=\lambda_N$.
We check Lyapunov's conditions to conclude (\ref{CLT}) and (\ref{PaN}).

\textit{Case} 2: $1-\min(2 \beta,\frac{4 \beta}{3-\tau}) < \zeta\leq
1-\beta$,
$\tau< \frac{\beta}{1-\zeta-\beta}$,
$\mu_N = (x-y) \sqrt{2 \log N}$
where $x=\sqrt{1-\zeta}$ and $y=\sqrt{(1+\tau)(1-\zeta-\beta)}$. Let
\begin{eqnarray*}
{\cal N}^0 & = & \bigl\{ n\dvtx I_n=0 \mbox{ or }
|Y_{n1}| < x \sqrt{2 \log N} \bigr\} ,
\\
{\cal N}^1 & = &\bigl \{
n\dvtx I_n=1 \mbox{ and } |Y_{n1}| \geq x \sqrt{2 \log N}
\bigr\},
\\
L_1^h & = & \prod_{n \in{\cal N}^h}
L_{n1},\qquad h=0,1.
\end{eqnarray*}
Check that $m_0 := E_1(L_{n1}|I_n=0) = 1$ and
%
%e7.2 #&#
\begin{eqnarray}
\label{m01} m_1 & := & E_1 \bigl[1+(L_{n1}-1)
{\mathbf I}_{\{ |Y_{n1}| \leq x \sqrt{2 \log N}
\}}|I_n=1 \bigr]
\nonumber
\\[-8pt]
\\[-8pt]
\nonumber
& = & 1+ \bigl[C_1+o(1) \bigr] \pi_N
N^{\zeta-1+2 \beta}/\sqrt{\log N}
\end{eqnarray}
for some $C_1>0$, hence
%
%e7.3 #&#
\begin{eqnarray}
\label{62a} E_1 \bigl(L_1^0 \bigr) & = &
\bigl[(1-\pi_N) m_0 +\pi_N m_1
\bigr]^N\nonumber
\\
& = & \bigl\{ 1+ \bigl[C_1+o(1) \bigr]
\pi_N^2 N^{\zeta-1+2 \beta}/\sqrt{\log N} \bigr\}
^N
\\
\nonumber
& = & \exp \bigl\{ \bigl[C_1+o(1) \bigr] N^{\zeta-2 \varepsilon}/
\sqrt{\log N} \bigr\}.
\end{eqnarray}

Next, we apply $\max_{n \in{\cal N}^1} |Y_{n1}| = O_p(\sqrt{\log N})$
to show that
%
%e7.4 #&#
\begin{eqnarray}
\label{62b} \log L_1^1 &= &O_p
\bigl( \bigl(\# {\cal N}^1 \bigr) \sqrt{\log N} \bigr)
\nonumber
\\[-8pt]
\\[-8pt]
\nonumber
& =&  O_p \bigl( N \pi_N \Phi \bigl(-
\sqrt{2(1-\zeta-\beta) \log N} \bigr) \sqrt{\log N} \bigr) = O_p
\bigl(N^{\zeta-\varepsilon} \bigr).
\end{eqnarray}
By (\ref{62a}), (\ref{62b}) and $L_1 = L_1^0 L_1^1$, (\ref{L1}) holds
for $\lambda_N= \exp(N^{\zeta-\varepsilon} \log N)$.
For $i \geq2$, let
\[
\wtd L_{ni} = L_{ni} {\mathbf I}_{\{ |Y_{ni}| \leq x \sqrt{2 \log N} \}},\qquad
\wtd L_i = \prod_{n=1}^N
\wtd L_{ni}.
\]
Then by a change-of-measure argument,
\begin{eqnarray*}
E_0(\wtd L_{ni}) & = & 1 + \pi_N
\bigl[P_1 \bigl\{ |Y_{n1}| \leq x \sqrt{2 \log N}
\bigr\}-P_0 \bigl\{ |Y_{n1}| \leq x \sqrt{2 \log N} \bigr\} \bigr]
\\
& =
& 1- \bigl[C_2+o(1) \bigr] N^{-\zeta-1-\varepsilon}/\sqrt{\log N},
\end{eqnarray*}
where $C_2 = \sqrt{1+\tau}/(2y \sqrt{\pi})$, therefore
\[
\kappa_N := E_0(\wtd L_i) = \exp \bigl\{-
\bigl[C_2+o(1) \bigr] N^{\zeta-\varepsilon
}/\sqrt {\log N} \bigr\}.
\]
Moreover, $E_0(\wtd L_{ni}^2) \geq[E_0(\wtd L_{ni})]^2$ for $i \geq
2$, and therefore, by (\ref{vNVar}),
\[
\operatorname{Var}_0(\wtd L_i) \geq
\bigl[1+o_p(1) \bigr] \kappa_N^2
\operatorname{Var}_0(\wtd L_{1i}).
\]
By (\ref{m01}) and a change-of-measure argument,
\[
\operatorname{Var}_0(\wtd L_{1i}) \sim E_0
\bigl(\wtd L_{1i}^2 \bigr) = E_1(\wtd
L_{11}) \sim C_1 N^{\zeta-1-2 \varepsilon}/\sqrt{\log N}.
\]
We check Lyapunov's condition to conclude (\ref{Lya2}) and (\ref{PaN}).

%s7.1 #&#
\subsection{Optimal detection using the penalized BJ test}\label{sec7.1}

By Lemma~\ref{lem1}, setting $h_N = 2 \log N$ leads to $P(\mbox{Type I
error}) \rightarrow0$.
To show $P(\mbox{Type II error}) \rightarrow0$, it suffices to find
$(j^*,j^*+\ell^*) \in B_{r,T}$ such that
$j_T \leq j^* < j^*+\ell^* \leq j_T+\ell_T$, $1-\ell^*/\ell_T =
O(r^{-1/2})$ and
%
%e7.5 #&#
\begin{equation}
\label{K8} P_1 \bigl\{ K(\bar S_N,t_N)
\geq N^{\zeta-1+\delta} \bigr\} \rightarrow1
\end{equation}
for some $0 < t_N < 1$ and $\delta> 0$, where
$\bar S_N = N^{-1} \sum_{n=1}^N {\mathbf I}_{\{ Y_{n \ell^* j^*} \geq
z(t_N) \}}$ and $\pi_N = N^{-\beta+\varepsilon}$,
$0 < \varepsilon< \beta$.

\textit{Case} 1(a): $0 \leq\zeta\leq1-2 \beta$, $\mu_N=0$. Let $t_N = \Phi
(-N^{-\varepsilon/2})$. Then
\[
E_1 \bar S_N - t_N = \pi_N
\bigl[ \Phi \bigl(-N^{-\varepsilon/2}/\sqrt{1+\tau } \bigr)-\Phi
\bigl(-N^{-\varepsilon/2} \bigr) \bigr] \sim C_3 N^{-\beta+\varepsilon/2},
\]
where $C_3=\frac{1}{\sqrt{2 \pi}}(1-\frac{1}{\sqrt{1+\tau}})$,
and since
\[
\operatorname{Var}_1(\bar S_N) \rightarrow(4N)^{-1}
= o \bigl((E_1 \bar S_N-t_N)^2
\bigr),
\]
therefore $\bar S_N-t_N \stackrel{p}{\rightarrow} C_3 N^{-\beta
+\varepsilon
/2}$. Since
$K(t,x) \sim2(t-x)^2$ when $t \rightarrow\frac{1}{2}$ and $x
\rightarrow\frac{1}{2}$,
(\ref{K8}) holds for $\delta=\varepsilon/2$.

\textit{Case} 1(b): $1-2 \beta< \zeta< 1-\frac{4 \beta}{3-\tau}$
($\Rightarrow
\tau<1$), $\mu_N
= C_4 \sqrt{\log N}$, where $C_4 = \sqrt{(1-\tau)(2 \beta+\zeta-1)}$.
Let $t_N = \Phi(-2 \mu_N/(1-\tau))$. Then
%
%e7.6 #&#
\begin{equation}
\label{8tN} t_N \sim C_5 N^{-2C_4^2/(1-\tau)^2}/\sqrt{
\log N},
\end{equation}
where $C_5 = (1-\tau)/(C_4 \sqrt{8 \pi})$, and
%
%e7.7 #&#
\begin{eqnarray}
\label{8diff} E_1 \bar S_N -t_N & = &
\pi_N \biggl[ \Phi \biggl( -\frac{\mu_N \sqrt
{1+\tau
}}{1-\tau} \biggr)
-t_N \biggr]
\nonumber
\\[-8pt]
\\[-8pt]
\nonumber
& \sim& C_6 N^{-C_4^2(1+\tau)/[2(1-\tau)]^2-\beta+\varepsilon}/\sqrt {\log N},
\end{eqnarray}
where $C_6 = (1-\tau)/(C_4 \sqrt{2 \pi(1+\tau)})$.

We claim that
%
%e7.8 #&#
\begin{equation}
\label{8claim} \operatorname{Var}_1 \bar S_N \biggl(
\sim N^{-1} \biggl[ t_N + \pi_N \Phi \biggl( -
\frac{\mu_N \sqrt{1+\tau}}{1-\tau} \biggr) \biggr] \biggr) = o \bigl((E_1 \bar
S_N)^2 \bigr).
\end{equation}
By (\ref{8tN})--(\ref{8claim}) and $t_N = o(E_1 \bar S_N)$,
%
%e7.9 #&#
\begin{eqnarray}
\label{8quad} (\bar S_N-t_N)^2/(2
t_N) & \stackrel{p} {\sim} & C_7 N^{C_4^2/(1-\tau)-2
\beta+2 \varepsilon}/
\sqrt{\log N}
\nonumber
\\[-8pt]
\\[-8pt]
\nonumber
& = & C_7 N^{\zeta-1+2 \varepsilon}/\sqrt{\log N}
\end{eqnarray}
for some $C_7>0$. Check that the inequality $-\frac{2C_4^2}{(1-\tau)^2}
> \zeta-1$ reduces to $\zeta< 1-
\frac{4 \beta}{3-\tau}$. Therefore by (\ref{8tN}), $t_N \sim C_5
N^{\zeta-1+2 \delta}/\sqrt{\log N}$
for some $\delta> 0$, and the root of $y^2/(2 t_N) = N^{\zeta
-1+\delta
}$ satisfies $y=o(t_N)$.
Since
$K(x,y) \sim\frac{(x-t)^2}{2t}$ as $t \rightarrow0$ and $\frac{x}{t}
\rightarrow1$,
(\ref{8quad}) implies (\ref{K8}).

It remains to show (\ref{8claim}) by comparing the leading exponent in
$N$ of the terms.
That is, it remains to show that
\[
-1+ \max \biggl( -\frac{2C_4^2}{(1-\tau)^2}, -\beta+\varepsilon-\frac
{C^2(1+\tau)}{2(1-\tau)^2}
\biggr) < -\frac{C_4^2(1+\tau)}{(1-\tau)^2}-2 \beta+2 \varepsilon,
\]
summarized as $-1+\max(A,B)< D$. The inequality
$-1+A<D$ reduces to $\zeta> -2 \varepsilon$, which holds trivially,
whereas $-1+B<D$ reduces to
\[
\zeta< \frac{3-\tau}{1+\tau} \biggl(1-\frac{4 \beta}{3-\tau} \biggr) +
\frac{2(1-\tau)}{1+\tau} \varepsilon,
\]
which holds because $\frac{3-\tau}{1+\tau} > 1$ when $\tau< 1$,
and it is assumed that $\zeta< 1-\frac{4 \beta}{3-\tau}$.
Hence (\ref{8claim}) holds.

\textit{Case} 2: $1-\min(2 \beta, \frac{4 \beta}{3-\tau}) \leq\zeta<
1-\beta$,
$\tau< \frac{\beta}{1-\zeta-\beta}$,
$\mu_N = (x-y) \sqrt{2 \log N}$ where $x=\sqrt{1-\zeta}$ and
$y=\sqrt
{(1+\tau)(1-\zeta-\beta)}$.
Let $t_N = \Phi(-x \sqrt{2 \log N}) [\sim
N^{\zeta-1}/(2x \sqrt{\pi\log N})]$. Then
\begin{eqnarray*}
E_1 \bar S_N & = & (1-\pi_N)
t_N+\pi_N \Phi \bigl(-\sqrt{2(1-\zeta -\beta ) \log N}
\bigr)
\\
& \sim& C_8 N^{\zeta-1+\varepsilon}/\sqrt{\log N}\qquad (\gg
t_N),
\end{eqnarray*}
where $C_8 = (2 \sqrt{\pi(1-\zeta-\beta)})^{-1}$. Moreover,
\[
\mathrm{ Var}_1 \bar S_N \sim E_1 \bar
S_N/N \sim C_8 N^{\zeta-2+\varepsilon
}/\sqrt{\log N} = o
\bigl((E_1 \bar S_N)^2 \bigr).
\]
Therefore $\bar S_N \stackrel{p}{\sim} C_8 N^{\zeta-1+\varepsilon
}/\sqrt
{\log N}$, and since
$K(x,t) \sim x \log(\frac{x}{t})$ as $x \rightarrow0$ and $\frac{x}{t}
\rightarrow\infty$,
%
%e7.10 #&#
\begin{equation}
\label{C9} K(\bar S_N, t_N) \stackrel{p} {\sim}
\bar S_N \log(\bar S_N/t_N) \stackrel{p} {
\sim} C_9 N^{\zeta-1+\varepsilon} \sqrt{\log N}
\end{equation}
for some $C_9 > 0$, (\ref{K8}) holds for $\delta=\varepsilon$. By similar
arguments, (\ref{K8}) holds for $\delta< \varepsilon$ when
$t_N \sim N^{\zeta-1}$.

\textit{Case} 3: $1-\min(2 \beta,\frac{4 \beta}{3-\tau}) \leq\zeta<
1-\beta$,
$\tau\geq\frac{\beta}{1-\zeta-\beta}$,
$\mu_N=0$. Let $t_N= \Phi(- \sqrt{2(1-\zeta) \log N})\ [\sim
N^{\zeta
-1}/(2x \sqrt{\pi\log N})]$. Then
%
%e7.11 #&#
\begin{equation}
\label{3SN} E_1 \bar S_N \sim\pi_N \Phi
\biggl(-\sqrt{\frac{2(1-\zeta) \log
N}{1+\tau
}} \biggr) \sim C_{10} N^{-\beta+\varepsilon
-(1-\zeta)/(1+\tau)}
\end{equation}
for some $C_{10}>0$. Since $\tau\geq\frac{\beta}{1-\zeta-\beta}$,
therefore the exponent of $N$ in
(\ref{3SN}) is at least $\zeta-1+\varepsilon$, and so $E_1 \bar S_N \gg t_N$.
We apply the first relation in (\ref{C9}) to conclude~(\ref{K8}),
for both $t_N = \Phi(-\sqrt{2(1-\zeta) \log N})$ and $t_N \sim
N^{\zeta-1}$.

%s7.2 #&#
\subsection{Optimal detection using the penalized HC test}\label{sec7.2}

By Lemma~\ref{lem2}, setting $h_N=2 \log N$ leads to $P(\mbox{Type I
error}) \rightarrow0$.
Let
\[
t_N = \cases{ %
\Phi
\bigl(-N^{-\varepsilon/2} \bigr), &\quad $\mbox{if } \zeta\leq1-2 \beta,$
\vspace*{2pt}\cr
\displaystyle \Phi
\biggl(-2 \sqrt{\frac{2 \beta-1+\zeta}{1-\tau} \log N} \biggr), & \quad
$\mbox{if } \displaystyle 1-2 \beta< \zeta
\leq1- \frac{4 \beta}{3-\tau},$
\vspace*{2pt}\cr
\displaystyle\frac{s_{\ell^* T}}{N}\ \bigl(\sim N^{\zeta-1}
\bigr), &\quad
$\mbox{if }\displaystyle 1-\min \biggl(2 \beta, \frac{4 \beta}{3-\tau} \biggr) < \zeta
\leq1-\beta.$}
\]
It was shown in (\ref{K8}) that in each case above,
$P_1 \{ K(\bar S_N,t_N) \geq N^{\zeta-1+\delta} \} \rightarrow1$ for
some $\delta>0$.
Since $\frac{x-t}{\sqrt{t(1-t)}} \geq\sqrt{2K(x,t)}$ for $x \geq t$
and $h_N+\sqrt{s_{\ell^* T} \log s_{\ell^* T}} =
o(N^{(\zeta+\delta)/2})$, $P(\mbox{Type II error}) \rightarrow0$.

%s7.3 #&#
\subsection{Proof of Corollary \texorpdfstring{\protect\ref{cor2}}{2}}\label{sec7.3}

Consider first the penalized HC test.
By Lem\-ma~\ref{lem2}, setting $h_N=2 \log N$ leads to $P(\mbox{Type I
error}) \rightarrow0$.
In the case $\pi_N = N^{-\beta+\varepsilon}$, the arguments above and in
Theorem~\ref{thmm2a} show that
\[
P_1 \biggl\{ \frac{\bar S_N-t_N}{\sqrt{t_N(1-t_N)/N}} \geq\sqrt {2N^\delta\log
\bigl(T/\ell^* \bigr)} \biggr\} \rightarrow1 \qquad \mbox{for some } \delta> 0.
\]
By (\ref{poly2}), $h_N+\sqrt{s_{\ell^* T} \log s_{\ell^* T}} =
o(\sqrt
{N^\delta})$, and
therefore $P(\mbox{Type II error}) \rightarrow0$, and (\ref{PP0})
holds. By similar arguments, (\ref{PP0}) holds for the penalized
BJ test.

%sA #&#
\begin{appendix}\label{app}
%sA #&#
\section{Verification of Lyapunov's condition}\label{appA}

We check in particular Lyapunov's condition to conclude (\ref{CLT}).
Let $\delta>0$ to be specified. It
follows from Taylor's expansion that
%
%eA.1 #&#
\begin{equation}
\label{A1} (1+u)^{2 +\delta} \leq1+(2+\delta)u+C^* u^2,\qquad
 |u| \leq1/2,
\end{equation}
for some $C^*>0$ chosen large enough. If $u>\frac{1}{2}$, then
%
%eA.2 #&#
\begin{equation}
\label{A2} (1+u)^{2 + \delta} \leq \biggl[ \sup_{v > 1/2}
\biggl( \frac{1+v}{v} \biggr)^{2+\delta} \biggr] u^{2+\delta}
=(3u)^{2+\delta}.
\end{equation}
By combining (\ref{A1}) and (\ref{A2}), we conclude that
%
%eA.3 #&#
\begin{equation}
\label{A3} (1+u)^{2+\delta} \leq1+(2+\delta)u + C^* u^2 +
C_2 |u|^{2+\delta} \qquad\mbox {for all } u \geq-1/2,
\end{equation}
where $C_2=3^{2+\delta}$. We apply (\ref{A3}) with $u=\pi_N[\exp
(\mu_N
Y_{n2}-\mu_N^2/2)-1]$ on (\ref{LNi})
to show that
%
%eA.4 #&#
\begin{eqnarray}
\label{Edelta}&& E_0 \bigl(L_2^{2+\delta} \bigr)
\nonumber
\\[-8pt]
\\[-8pt]
\nonumber
&&\qquad\leq \bigl[1+C^* p_N^2 \exp \bigl(\mu_N^2
\bigr)+C_2 \pi _N^{2+\delta
} \exp \bigl(
\mu_N^2 (1+\delta) (2+\delta)/2 \bigr)
\bigr]^N.
\end{eqnarray}
Since
$\mu_N = \sqrt{\log(1+N^{2 \beta-1+\zeta})}$ and $\pi_N =
N^{-\beta
-\varepsilon}$, by (\ref{Edelta}),
%
%eA.5 #&#
\begin{equation}
\label{A4} E_0 \bigl(L_2^{2+\delta} \bigr)
\leq \exp \bigl\{ \bigl[1+o(1) \bigr] \bigl(C^* N^{\zeta-2 \varepsilon}+C_2
N^{\zeta-2 \varepsilon+\kappa} \bigr) \bigr\},
\end{equation}
where $\kappa=-\delta(\beta+\varepsilon)+\frac{3 \delta+ \delta^2}{2}(2
\beta-1+\zeta)$. Let $\delta>0$ be
small enough such that $\kappa< \varepsilon$. Since $i_T-1 \sim\exp
(N^{\zeta}-1)$
and $L_i-1 \geq-1$, we get from (\ref{A4}) that
\[
(i_T-1) E_0 (L_i-1)^{2+\delta} \leq
\exp \bigl\{ \bigl[1+o(1) \bigr] \bigl(N^\zeta+C_2
N^{\zeta-
\varepsilon} \bigr) \bigr\}.
\]
On the other hand, $\operatorname{Var}_0(L_i)=\lambda_N-1$ and
\[
\bigl[(i_T-1) (\lambda_N -1) \bigr]^{1+\delta/2} =
\exp \bigl\{ \bigl[1+\delta/2+o(1) \bigr] \bigl(N^\zeta +
N^{\zeta-2 \varepsilon} \bigr) \bigr\},
\]
so Lyapunov's condition is satisfied.

%sA #&#
\section{Quadratic bounds for the function $K$}\label{appB}

For given $\ell$, $T$, $N$, $\frac{s_{\ell T}}{N} \leq t \leq
\frac{1}{2}$ and $0 < \gamma\leq\frac{c_{N \ell}}{N}$
(recall that $c_{N \ell} = \log N + s_{\ell T} + 3 \log s_{\ell T}$),
let $x_t > t$ be such that
%
%eA.1 #&#
\begin{equation}
\label{Qxt} Q(x_t,t)\ \biggl[ :=\frac{(x_t-t)^2}{2t(1-t)} \biggr] =
\gamma.
\end{equation}

We claim that
%
%eA.2 #&#
\begin{equation}
\label{KQ2} K(x_t,t) \leq Q(x_t,t) \leq \Bigl(3
\sqrt{4+ s_{\ell T}^{-1} \log N} \Bigr)
K(x_t,t).
\end{equation}
The left inequality of (\ref{KQ2}) is known. To obtain the
right inequality, first note that by (\ref{Qxt}),
\[
x_t = t(1+y)\qquad \mbox{where } y=\sqrt{2 \gamma(1-t)/t}.
\]
Since $\frac{d}{dx}[(1-x) \log(\frac{1-x}{1-t})] = -1-\log(
\frac{1-x}{1-t}) \geq-1$ for $x \geq t$, therefore
%
%eA.3 #&#
\begin{equation}
\label{C2} K(x_t,t) \geq x_t
\log(x_t/t)-(x_t-t) = t \bigl[f(y) \bigr],
\end{equation}
where $f(y) = (1+y) \log(1+y)-y$.
Check that $\frac{d}{dy} f(y)=\log(1+y)$, and that $\frac{d^2}{dy^2}
f(y)=\frac{1}{1+y}$.

For $\gamma\leq t$, apply (\ref{Qxt}), (\ref{C2}) and
$f(y) \geq\frac{y^2}{4}$ on $0 < y \leq1$ to show that
\[
\frac{Q(x_t,t)}{K(x_t,t)} \leq\frac{4 \gamma}{ty^2} = \frac{2}{1-t} \leq4.
\]
The right inequality of (\ref{KQ2}) holds.

For $\gamma> t$, apply (\ref{Qxt}), (\ref{C2}) and $f(y) \geq y/3$
on $y>1$ to show that
\[
\frac{Q(x_t,t)}{K(x_t,t)} \leq\frac{3 \gamma}{t y} = 3 \sqrt{\frac{\gamma}{2t(1-t)}} \leq3
\sqrt{\frac{\gamma}{t}} \leq3 \sqrt{\frac{c_{N \ell}}{s_{\ell T}}}.
\]
The right inequality of (\ref{KQ2}) again holds.

%sA #&#
\section{Detectability of nonaligned signals}\label{appC}

Let $X_{nt}=\mu_{nt}+Z_{nt}$, where $Z_{nt}$ are i.i.d. $\mathrm{N}(0,1)$.
Assume that there exists $1 \leq\ell_T \leq T$ such that for the $n$th
sequence,
$1 \leq n \leq N$,
there is an unknown interval $(j_{nT}, j_{nT}+\ell_T]$ with probability
$\pi_N > 0$ of having an elevated mean
%
%eA.1 #&#
\begin{eqnarray}
\label{D.0} \mu_{nt} & = & \cases{ %
\mu_N I_n/\sqrt{\ell_T}, & \quad$\mbox{if }
j_{nT} < t \leq j_{nT}+\ell_T,$
\vspace*{2pt}\cr
0, & \quad$\mbox{otherwise},$}
\nonumber
\\[-8pt]
\\[-8pt]
\nonumber
I_n & \sim&  \operatorname{Bernoulli}(\pi_N),
\end{eqnarray}
with $\mu_N > 0$ and the $I_n$'s and $Z_{nt}$'s jointly independent.
This model
is distinct from (\ref{munt}) in that we do not now assume that the
signals are aligned.

We claim that if $\pi_N = N^{-\beta-\varepsilon}$ for some $0 < \beta< 1$
and $\varepsilon> 0$,
%
%eA.2 #&#
\begin{equation}
\label{TlT} T/\ell_T \sim N^\zeta\qquad\mbox{for some }
\zeta> \max(0,1-2 \beta),
\end{equation}
and $\mu_N = \sqrt{(2 \log N)(\zeta+1) \rho^*(\frac{\zeta+\beta
}{\zeta
+1})}$, then there is no test that can achieve at all
such $j_{nT}$,
%
%eA.3 #&#
\begin{equation}
\label{D2} P(\mbox{Type I error}) + P(\mbox{Type II error}) \rightarrow0.
\end{equation}
Note that though $\mu_N \,\dot{\sim}\, \sqrt{\log N}$, as in the boundary
for cases 1(b) and 2 in (\ref{bNbz}), the growth of
$T/\ell_T$ that is allowed in (\ref{TlT}) is considerably smaller than
that of (\ref{C1}).

As in the proof of Theorem~\ref{thmm1}, set $i_T = \lfloor T/\ell_T
\rfloor-1$, so that $i_T \sim N^\zeta$. Let $Y_{n1} = \ell_T^{-1/2}
(X_{n,j_T+1} + \cdots+ X_{n,j_T+\ell_T})$ and each
$Y_{ni}$, $2 \leq i \leq i_T$ be of the form $\ell
_T^{-1/2}(X_{n,j+1}+\cdots+X_{n,j+\ell_T})$, with all $(j,j+\ell_T]$
disjoint from each other, and from $(j_{nT},j_{nT}+\ell_T]$. Assume
without loss of generality each $j_{nT}$ is equally likely to
take one of the $i_T$ possible values spaced at least $\ell_T$ apart,
as given above.
Then when $\varepsilon=0$, the detection of
(nonaligned) signals satisfying (\ref{TlT}) is at least as difficult
as detecting a mixture of $\sim\!N^{\zeta+1}$
normal random variables $\{ Y_{ni}\dvtx 1 \leq i \leq i_T \}$, with a
sparse fraction
$\sim\! N^{-(\zeta+\beta)}\ [=
(N^{\zeta+1})^{-{(\zeta+\beta)}/{(\zeta+1)}}]$ of them
having mean $\mu_N$. Therefore by the results in \cite{Ing97,Ing98},
the critical detectable $\mu_N$ is
$\sqrt{(2 \log N^{\zeta+1})\rho^*(\frac{\zeta+\beta}{\zeta
+1})}$. Hence
when $\varepsilon> 0$, (\ref{D2}) cannot be achieved.
Note that the assumption $\zeta> 1-2 \beta$ in (\ref{TlT})
implies that $\frac{\zeta+\beta}{\zeta+1}
> \frac{1}{2}$, so $\rho^*(\frac{\zeta+\beta}{\zeta+1})$ is well-defined.
\end{appendix}

\section*{Acknowledgments} We would like to thank an Associate Editor
and two referees for their comments
that have led to a more realistic model setting and the adaptive
optimal test statistics in this paper.
% imsref loaded by akundreckaite, 2015-03-24 13:25:31
% imsref loaded by akundreckaite, 2015-03-24 14:37:09

%\begin{appendix}
%\section{}
%\end{appendix}

%\section*{Acknowledgments}

%\begin{supplement}[id=suppA]
%\slink[doi]{10.1214/00-AOSXXXXSUPP} %[doi,text={...}] - jei reikia
%suskaldyti doi
%\sfilename{aosXXXX\_supp.pdf}
%\end{supplement}

%\begin{thebibliography}{99}
%\bibitem[\protect\citeauthoryear{}{}]{r1}
%\bibitem{r1}
%\end{thebibliography}

\printaddresses
\end{document}